\documentclass[final, 3p, times, review]{elsarticle}
\usepackage[utf8]{inputenc}
\usepackage{geometry}
\usepackage{framed,multirow,hyperref}
\usepackage{geometry}
\usepackage{amssymb}
\usepackage{amsthm}
\usepackage{amsmath}
\usepackage{amsfonts,amssymb}
\usepackage{enumerate}
\usepackage{booktabs}

\usepackage{graphicx}
\usepackage{multirow}
\usepackage{mathtools}
\usepackage{color}
\usepackage{algorithm}
\usepackage{algpseudocode}
\usepackage{appendix}
\usepackage{stmaryrd}
\usepackage{subfigure}
\usepackage{mathrsfs}
\usepackage{threeparttable}
\usepackage{diagbox}
\usepackage{hyperref}
\usepackage{makecell}

\newtheorem{proposition}{Proposition}[section]

\newtheorem{remark}{Remark}[section]

\newcommand{\R}{\mathbb{R}}

\newcommand{\norm}[1]{\left\lVert#1\right\rVert} 
\renewcommand{\vec}[1]{\boldsymbol{#1}}

\renewcommand{\vec}[1]{\boldsymbol{#1}}

\usepackage{ulem}

\begin{document}

\begin{frontmatter}
\title{Fine-Tuning DeepONets to Enhance Physics-informed Neural Networks for solving Partial Differential Equations}
\author[1]{Sidi Wu\corref{cor}}
\address[1]{School of Mathematical Sciences, Peking University, Beijing 100871, China} 
\cortext[cor]{Corresponding author. Email addresses: wsd@pku.edu.cn}
\date{}

\begin{abstract} 
Physics-Informed Neural Networks (PINNs) have emerged as powerful tools for solving partial differential equations (PDEs). However, training PINNs from scratch is often computationally intensive and time-consuming. To address this problem, we propose a parameter-efficient approach that fine-tunes pre-trained DeepONet models within the PINN framework (FTO-PINN), enabling more efficient meshless PDE solving. Specifically, we freeze the weights of the pre-trained DeepONet model and fine-tune the output of the branch net by incorporating a small number of new trainable parameters, which can be quickly determined using least-squares techniques. Additionally, we introduce trunk net expansions and low-rank adaptation strategies to further enhance the performance of FTO-PINN. The effectiveness of our proposed method is demonstrated through a series of numerical experiments across various types of PDEs. FTO-PINN significantly reduces the training time of vanilla PINNs while maintaining comparable accuracy, and outperforms DeepONet, which is pre-trained on general function data, in both fidelity and generalization capabilities.

\end{abstract}
\end{frontmatter}

\section{Introduction}
 
The remarkable success of deep learning in diverse fields has led to its increasing adoption in scientific computing. A particularly promising avenue is the use of neural networks to solve partial differential equations (PDEs) without relying on traditional mesh-based methods. Among these approaches, Physics-Informed Neural Networks (PINNs) \cite{raissi2019physics,karniadakis2021physics} have emerged as a powerful paradigm. PINNs leverage neural networks to approximate PDE solutions and train them by minimizing a loss function that combines data fitting and physical equation enforcement.
Recent advancements in PINNs have demonstrated their potential in various scientific domains, including fluid mechanics \cite{raissi2020hidden}, systems biology \cite{yazdani2020systems}, biomedicine \cite{wu2022interfaced,ying2024an}, and even inverse problems of stochastic PDEs \cite{zhang2019quantifying}.  Meanwhile, robust general-purpose algorithms have become available \cite{lu2021deepxde} and our theoretical understanding of these methods is also being developed \cite{shin2024chapter, lu2022machine, gao2023gradient}.

Despite the aforementioned success, solving PDEs with PINNs involves optimizing non-convex objective functions, resulting in computationally intensive and time-consuming training processes. This issue is prevalent across various deep learning-based methods. For instance, the recent language model GPT-3, with its 175 billion parameters, required 355 GPU years to train \cite{gpt2020report}. Transfer learning, a technique widely used in computer vision and natural language processing, offers a solution to mitigate the computational challenges. By fine-tuning pre-trained models trained on general-domain data, it enables efficient adaptation to specific tasks without starting the training process from scratch. Neural network-based PDE solvers have also benefited from transfer learning \cite{weinan2018deep,goswami2020transfer,chakraborty2021transfer,chen2021transfer,niaki2021physics,xu2023transfer}, enabling the reuse of knowledge from previously trained models. However, this task typically lacks a general-purpose dataset; each PDE typically requires its own specialized dataset, which hinders the direct application of fine-tuning techniques in transfer learning.

Recent advancements in neural networks have led to the development of operator learning methods capable of directly mapping between infinite-dimensional function spaces. Examples include Deep Operator Networks (DeepONet) \cite{lu2021learning} and Fourier Neural Operators \cite{li2021fourier}. Given their training on high-level functional data, these models hold significant potential for adapting to and solving specific PDE problems. In this study, we propose a novel method for fine-tuning pre-trained DeepONets within the PINN framework (FTO-PINN) to efficiently solve specific PDE problems. Specifically, we freeze the pre-trained model and introduce a new set of trainable parameters to fine-tune the output of the branch net. This strategy substantially enhances parameter efficiency while preserving the representational capacity of the original model. Instead of using propagation-type methods such as Adam, the proposed FTO-PINN enables rapid determination of trainable weight values through least-squares techniques, significantly accelerating the training process. 
Additionally, we introduce strategies to expand the trunk net and incorporate trainable matrices into linear layers of trunk net to enhance the accuracy of FTO-PINN, ensuring high fidelity even for new PDEs beyond the domain for which the DeepONet was originally trained and performs adequately. Several key advantages of our proposed method are summarized below and will be documented in detail later.
\begin{itemize}
\item FTO-PINN demonstrates improved computational efficiency compared to standard PINNs, while achieving comparable or superior accuracy.
\item Compared to DeepONet pre-trained on general function data, FTO-PINN demonstrates superior fidelity and generalization.
\item FTO-PINN is parameter-efficient and applicable to a wide range of PDEs, including linear, nonlinear, and interface problems.
\end{itemize}

The rest of the paper is organized as follows: Section \ref{sec:related works} provides a brief overview of related work. Section \ref{sec:preliminary} presents the necessary preliminaries. Section \ref{sec:methodologies} details the FTO-PINN algorithm, and Section \ref{sec:numerical results} presents numerical experiments to validate the effectiveness of our proposed method.
Finally, Section \ref{sec:conclusion} concludes the paper.

\subsection{ Related works}\label{sec:related works}
Fine-tuning, within the realm of machine learning, denotes the process of taking a pre-trained model and further adjusting its certain parameters. It has become a fundamental transfer learning technique, particularly in the training of foundation models used for generative natural language processing and computer vision tasks, such as GPT-3 \cite{language2020tom} and CLIP \cite{radford2021learning}. Fine-tuning allows for the adaptation of a model pre-trained on general domains to a specific task, circumventing the need to start training from scratch \cite{ruder2018universal}. Variants of this approach include optimizing only a subset of the parameters \cite{beta2019jacob} or incorporating external modules to enhance downstream performance \cite{li2021prefix,lora2022hu}. This paper delves into scientific machine learning, focusing on fine-tuning of a high-level pre-trained model, specifically an operator approximator, to address a specific PDE problem.

The use of neural networks to solve PDEs dates back to the early 1990s \cite{psichogios1992hybrid, lagaris1998artificial}. In recent years, this approach has been further advanced, with notable methods including the deep Ritz method (DeepRitz) \cite{weinan2018deep}, the deep Galerkin method \cite{dgm2018dgm} and PINNs \cite{raissi2019physics}. These methods are particularly advantageous due to their meshless nature and ease of implementation, allowing for straightforward extension to accommodate irregular domains \cite{wu2022interfaced}, multiple domains \cite{jagtap2020conservative}, and high-dimensional PDEs \cite{han2018solving}. While the transition from numerical discretization schemes to mathematical optimization is promising, these methods need to be further improved in terms of computational cost and accuracy compared to traditional methods \cite{xu2023transfer,chen2022bridging}. One potential improvement is to use random-weighted models (RWMs) \cite{scardapane2017randomness,ying2024an}, also referred to as extreme learning machines \cite{huang2011extreme, dwivedi2020physics} or random feature methods \cite{chen2022bridging,sun2024local,xu2024subspace}. RWMs differ from general NN-based methods by simplifying the training process through the random initialization and fixation of hidden layer parameters, thereby focusing learning on the output layer. A key challenge in this line of work is the effective generation of random coefficients to enhance performance \cite{dong2021modified,zhang2024transferable}.

Neural operators leverage neural networks to directly learn the solution mappings between two infinite-dimensional function spaces. Notable examples include the Deep operator network (DeepONet) \cite{lu2021learning}, the Fourier neural operator \cite{li2021fourier}, and Meta-Mgnet \cite{chen2022meta}. These emerging methods offer innovative approaches for achieving rapid solutions to PDEs and provide new insights into physical phenomena. They have demonstrated strong performance across a range of applications, including multiscale bubble growth dynamics \cite{lin2021operator}, multiphysics and multiscale problems in hypersonics \cite{mao2021deepm} and multiphase flow \cite{wen2022u}. With the rapid growth of this field, many neural operators have been proposed to achieve higher performance or reduce the amount of data required for training \cite{wang2021learning,jin2022mionet,li2023fourier,wu2024solving,yang2023in}. Despite these successes, neural operators are trained on general function data, and may exhibit low accuracy in specific PDE-solving tasks, especially when applied to new PDEs beyond the scope of the training set \cite{xu2020neural}.

Several studies have investigated the fine-tuning techniques for neural network-based solutions. A common strategy is to initialize a network with the parameters of a model pre-trained on a related task and then fine-tune the entire network \cite{weinan2018deep,niaki2021physics,lippert2024transfer} or specific layers \cite{goswami2020transfer,chakraborty2021transfer,xu2023transfer} on the target task. This idea has also been extended to scenarios containing multiple pre-trained NNs \cite{wang2022mosaic}. The use of this warm-start initialization strategy can significantly reduce the overall training time of methods such as PINN and DeepRitz, especially when the solution of the target problem is close to that of the original problem where the model was trained. In addition to fine-tuning models trained on specific PDE problems, there are also works focusing on fine-tuning neural operators, such as DeepONet. For instance, \cite{goswami2022deep} constructs a new hybrid loss function to fine-tune the pre-trained DeepONet from one computational domain to another.
By leveraging underlying physical equations or new observations, \cite{lin2021operator, zhu2023reliable} propose fine-tuning the pre-trained DeepONet to achieve accurate predictions when the input lies outside the training data space. These existing fine-tuning techniques for PDE solvers typically transfer models trained for one task to another of a similar scale. In contrast, FTO-PINN aims to adapt a DeepONet model, pre-trained on high-level function data, to efficiently solve a specific PDE as a downstream task, achieved through a customized fine-tuning scheme tailored for the task.

\section{Preliminary\label{sec:preliminary}}
Consider a general PDE \begin{equation}
 \begin{aligned}
  \mathcal{L}(u(\mathbf{x})) &= f(\mathbf{x}), \ \mathbf{x}\in\Omega,\\
  \mathcal{B}(u(\mathbf{x})) &= h(\mathbf{x}), \ \mathbf{x}\in \partial\Omega.
 \end{aligned}
 \label{eq: general specific pde}
\end{equation}
Here, $u$ denotes the solution of interest, $\mathcal{L}$ is a general differential operator, $\mathcal{B}$ is an operator defined on $\partial \Omega$, and $f$ and $h$ are the given right-hand term and boundary condition, respectively. When a time-dependent PDE is involved, the computational domain $\Omega$ includes the temporal domain, and the time variable $t$ can be treated as a special component of the spatial variable $\mathbf{x}\in \mathbb{R}^{d_0}$. In this way, the initial conditions are regarded as special Dirichlet boundary conditions.  In the following, we assume that the solution $u$ to PDE \eqref{eq: general specific pde} exists and has some degree of continuity, i.e., its differentiability is at least of the same order as that of the PDE.

\subsection{Feedforward neural network}
Feedforward neural network (FNN) is one of the most popular neural network architectures \cite{raissi2019physics,weinan2018deep}. Mathematically, an $L$-layer FNN $u_{\vec{\theta}}: \mathbb{R}^{d_0}\rightarrow \mathbb{R}^{d_L}$ is recursively defined as follows:
\begin{equation}
\begin{aligned}
 &\psi_0 = \mathbf{x},\\
 &\psi_l =\sigma(\mathbf{W}_l\psi_{l-1}+\mathbf{b}_l),\; 1\leq l \leq L-1,\\
 &\psi_L= \mathbf{W}_L\psi_{L-1} + \mathbf{b}_L.
\end{aligned} 
\label{eq: fully connection network}
\end{equation}
Here, $\vec{\theta}:=\{\mathbf{W}_l,\mathbf{b}_{l}\}_{l=1}^{L}$ denotes the set of all trainable parameters, where $\mathbf{W}_l\in \mathbb{R}^{d_l\times d_{l-1}}$ and $\mathbf{b}_l\in \mathbb{R}^{d_l}$ are the weight matrix and bias vector in the $l$-th layer, $\sigma:\mathbb{R} \rightarrow \mathbb{R}$ is a nonlinear activation function applied element-wise to a vector. Commonly used activation functions include $\texttt{Sigmoid}(z)=1/(1+e^{-z})$, $\texttt{Tanh}(z)=(e^z-e^{-z})/(e^z+e^{-z})$ and $\texttt{ReLU}(z)=\max\{0, z\}$. FNNs are widely used in the field of scientific computing due to their exceptional ability to approximate complex functions \cite{pinkus1999approximation,guhring2021approximation} and generalize to unseen data \cite{lecun2015deep}.

\subsection{Physics-informed neural networks}
The main idea of physics-informed neural networks (PINNs) is to train a deep neural network $u_{\vec{\theta}}(\mathbf{x})$ to approximate the solution $u(\mathbf{x})$ by minimizing the residuals of the PDEs and the boundary conditions. Given the training data $\mathcal{T}=\mathcal{T}_{D}\cup\mathcal{T}_{B}$, 
which comprises two sets, $\mathcal{T}_{D}=\{\mathbf{x}_j^d\}_{j=1}^{N_D}\subset\Omega$ and $\mathcal{T}_{B}=\{\mathbf{x}_j^b\}_{j=1}^{N_{B}}\subset\partial\Omega$, these represent scattered points in the domain and on the boundary, respectively. The loss function of the PINN for solving problem  \eqref{eq: general specific pde} has the following specific form:
\begin{equation}
\text{Loss}(\vec{\theta}, \mathcal{T}) = \frac{1}{N_D}\sum_{j=1}^{N_D}\rho_j^d\left| \mathcal{L}\left( u_{\vec{\theta}}(\mathbf{x}_j^d)\right) -f(\mathbf{x}_j^d)\right|^2
+
\frac{1}{N_{B}}\sum_{j=1}^{N_{B}}\rho_j^b \left| \mathcal{B}\left( u_{\vec{\theta}}(\mathbf{x}_j^b)\right) -h(\mathbf{x}_j^b)\right|^2,
\label{eq: pinn loss}
\end{equation}
where $\rho_j^d>0$ and $\rho_j^b>0$ are the penalty weights at different collocation points, typically set to 1 by default. Derivatives in the loss can be handled via automatic differentiation \cite{baydin2017automatic}. Since the loss is highly nonlinear and nonconvex with respect to $\vec{\theta}$, gradient-based optimizers such as Adam and LBFGS are usually used to minimize the loss function.

\subsection{Deep operator networks}
Deep operator network (DeepONet) \cite{lu2021learning} and its extensions \cite{wang2021learning,wu2024solving} have emerged as a powerful tool for solving parametric PDEs. They learn solution operators $\mathcal{G}$ that directly map PDE parameters $v$ (e.g., source terms, coefficients, and boundary conditions) to the corresponding PDE solutions $u$. Specifically, a conventional DeepONet $\mathcal{G}_{\vec{\theta}}$ comprises two FNNs, referred to as the branch net $\mathcal{N}_{b}$ and the trunk net $\mathcal{N}_{t}$. The prediction of an input function $v$ evaluated at a point $\mathbf{x}$ by $\mathcal{G}_{\vec{\theta}}$ can be expressed as:
\begin{equation}
 \mathcal{G}_{\vec{\theta}}(v)(\mathbf{x}) = \left\langle \mathcal{N}_{b}(\overrightarrow{v}) , \mathcal{N}_{t}(\mathbf{x}) \right\rangle,
 \label{eq: DeepONet}
\end{equation}
where $\mathcal{N}_{b}: \R^m\rightarrow \R^I$, $\mathcal{N}_{t}: \R^{d_0}\rightarrow \R^I$ and $\langle \cdot, \cdot \rangle$ denotes the dot product in $\R^I$. The vector $\overrightarrow{v}:=[v(\mathbf{y}_1),\cdots,v(\mathbf{y}_m)]$ serves as the input to $\mathcal{N}_{b}$, comprising the values of the input function $v$ at a set of fixed points $\{\mathbf{y}_i\}_{i=1}^m$. In particular, when the trunk net $\mathcal{N}_{t}$ is equipped with a smooth activation function such as \texttt{Tanh} and \texttt{Sigmoid}, $\mathcal{G}_{\vec{\theta}}$ \eqref{eq: DeepONet} is capable of delivering a continuous and differentiable approximation of the solution function $u$. This advancement has spurred research into the development of the Physics-Informed DeepONet (PI-DeepONet), as detailed in \cite{wang2021learning}.
In addition, in combination with the domain decomposition technique, IONet \cite{wu2024solving}, a variant of DeepONet, efficiently handles operators with discontinuous input and output functions. In the scenario of two subdomains (i.e., $\Omega=\Omega_1\cup\Omega_2$), IONet has the following network structure
\begin{equation}
\mathcal{G}_{\mathbf{\vec{\theta}}}(v)(\mathbf{x}) = \left\{\begin{aligned}
&\mathcal{G}_{\mathbf{\vec{\theta}}}^1(v)(\mathbf{x}),\quad  \mathbf{x} \in \Omega_1,\\
&\mathcal{G}_{\mathbf{\vec{\theta}}}^2(v)(\mathbf{x}),\quad  \mathbf{x} \in \Omega_2,\\
\end{aligned}\right.
\label{eq: structure of ionet}
\end{equation}
where $\mathcal{G}^j_{\mathbf{\vec{\theta}}}(v) (\mathbf{x})= \langle b(\overrightarrow{v}),\mathcal{N}^j_{t}(\mathbf{x})\rangle$ with $j=1,2$, $b(\overrightarrow{v}) = 
\mathcal{N}_{b_1}(\overrightarrow{v}\big|_{\Omega_1})\odot\mathcal{N}_{b_2}(\overrightarrow{v}\big|_{\Omega_2})$, and  $\odot$ represents the Hadamard product.

\section{Methodology\label{sec:methodologies}}

In this section, we present the methodology for fine-tuning pre-trained DeepONets within the PINN framework, referred to as FTO-PINN, to effectively address a specific PDE problem \eqref{eq: general specific pde}.

\subsection{Fine-tuning with a pre-trained DeepONet}
 
Suppose we have a DeepONet model, denoted as $\mathcal{G}_{\vec{\theta}}$, parameterized by $\vec{\theta}$ and pre-trained on a general function dataset $\{(v_s, u_s)\}_{s=1}^S$. Here, the input function $v_s$ represents the parameters of the PDE \eqref{eq: general specific pde}, such as source terms, coefficients, and boundary conditions, while $u_s$ corresponds to the latent solution. We then fine-tune $\mathcal{G}_{\vec{\theta}}$ for a downstream task of solving a specific PDE with fixed parameters.

In the architecture of DeepONet, the trunk net $\mathcal{N}_{t}$ is responsible for extracting continuous input coordinates $\mathbf{x}$ where the output function is evaluated, while branch net $\mathcal{N}_{b}$ encodes the input function $v$ into latent representations. In the framework of FTO-PINN, we solely employ the trunk net component for PDE solving. During training, the trunk net is frozen and does not receive gradient updates. The weights previously output by the branch net are replaced with a set of new trainable parameters $\vec{\alpha} =[\alpha_1,\ldots,\alpha_I]\in\mathbb{R}^I$. Denote $\mathcal{N}_{t}(\mathbf{x}):=[t_1(\mathbf{x}), t_2(\mathbf{x}), \ldots, t_I(\mathbf{x})]^T\in\mathbb{R}^I$ and $\mathcal{N}_{b}\left(\overrightarrow{v}\right) := [b_1\left(\overrightarrow{v}\right), \ldots, b_I\left(\overrightarrow{v}\right)]^T\in\mathbb{R}^I$ as the outputs of the trunk net and branch net, 
then for $\mathcal{G}_{\vec{\theta}}(v)(\mathbf{x})=\sum_{i=1}^Ib_i(\overrightarrow{v})t_i(\mathbf{x})$, our modified forward pass yields:
\begin{equation}
u_{\vec{\alpha}}(\mathbf{x}): = \sum_{i=1}^I \alpha_i  t_i(\mathbf{x}).
\label{eq: general form of fto-pinn approximate solution}
\end{equation}
Formally, this is equivalent to fine-tuning the output of the branch net without altering its internal parameters. In contrast to existing  fine-tuning methods \cite{lin2021operator,zhu2023reliable} that directly update the parameters of the model itself, FTO-PINN substantially reduces the number of trainable parameters, leading to potential enhancements in training efficiency.

\begin{figure}[htbp] 
	\centering
\includegraphics[width=1.0\textwidth]{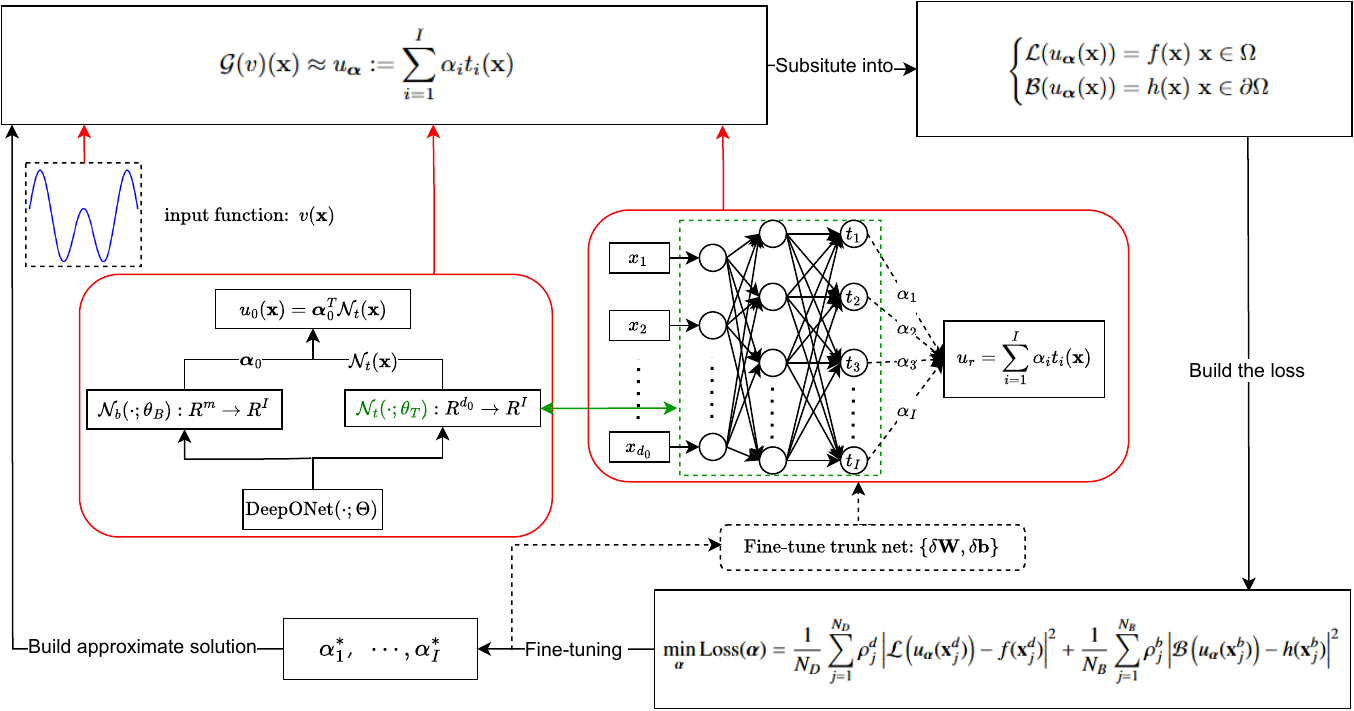}
	\caption{Schematic diagram of FTO-PINN learning procedure.}
\label{fig: framework}
\end{figure}

By inserting \eqref{eq: general form of fto-pinn approximate solution} into Eq. \eqref{eq: general specific pde} and following the principles of PINN, the task of identifying the parameters $\vec{\alpha}$ depends on the following optimization problem
\begin{equation}
\min_{\vec{\alpha}}\text{Loss}(\vec{\alpha})=\frac{1}{N_D}\sum_{j=1}^{N_D}\rho_j^d\left| \mathcal{L}\left( u_{\vec{\alpha}}(\mathbf{x}_j^d)\right)-f(\mathbf{x}_j^d)\right|^2
+
\frac{1}{N_{B}}\sum_{j=1}^{N_{B}}\rho_j^b \left| \mathcal{B}\left( u_{\vec{\alpha}}(\mathbf{x}_j^b)\right)-h(\mathbf{x}_j^b)\right|^2.
\label{eq: loss function of fto-pinn}
\end{equation}
The algorithm of our proposed FTO-PINN is illustrated in  Fig. \ref{fig: framework}. We next discuss how to determine the parameter $\vec{\alpha}$ for the specific PDEs \eqref{eq: general specific pde}.
\begin{remark} 
DeepONet has demonstrated effectiveness in learning nonlinear continuous operators and has produced numerous pre-trained DeepONet models for a wide range of PDE systems. For a specific target PDE, FTO-PINN can select a pre-trained DeepONet that has been trained to approximate solution operator of similar or related but not entirely identical PDEs.  
\end{remark}

\subsubsection{Linear PDEs\label{subsection: linear pde}}
Consider a linear PDE problem \eqref{eq: general specific pde}, where $\mathcal{L}$ and $\mathcal{B}$ are linear. Combining the framework of PINN, the goal of FTO-PINN is to find a set of $\vec{\alpha}=[\alpha_1,\cdots,\alpha_I]$ such that $u_{\vec{\alpha}}$ \eqref{eq: general form of fto-pinn approximate solution} minimizes the loss function \eqref{eq: loss function of fto-pinn}. By inserting $u_{\vec{\alpha}}$ into the loss \eqref{eq: loss function of fto-pinn}, we obtain a linear least-squares problem
\begin{equation*}
\min_{\vec{\alpha}}\text{Loss}(\vec{\alpha})=\frac{1}{N_D}\sum_{j=1}^{N_D}\rho_j^d\left| \sum_{i=1}^I \alpha_i \mathcal{L}\left( t_i(\mathbf{x}_j^d)\right) -f(\mathbf{x}_j^d)\right|^2
+
\frac{1}{N_{B}}\sum_{j=1}^{N_{B}}\rho_j^b \left|\sum_{i=1}^I\alpha_i \mathcal{B}\left( t_i(\mathbf{x}_j^b)\right)-h(\mathbf{x}_j^b)\right|^2,
\end{equation*}
which can be rearranged as follows:
\begin{equation}
\min_{\vec{\alpha}}\text{Loss}(\vec{\alpha}) = \|\vec{\rho}\mathbf{A} \vec{\alpha} - \vec{\rho} \mathbf{b}\|_2^2,
\label{eq: lstq system}
\end{equation}
where 
\begin{equation*}
 \mathbf{A} = 
 \left[
 \begin{array}{c}
  \mathbf{A}_1 \\
  \mathbf{A}_2 \\ 
 \end{array}
\right], \
\mathbf{A}_1= \left[
 \begin{array}{ccc}
  \mathcal{L}\left( t_1(\mathbf{x}_1^d)\right) &\cdots & \mathcal{L}\left( t_I(\mathbf{x}_1^d)\right)\\
  \vdots & \ddots & \vdots\\
   \mathcal{L}\left( t_1(\mathbf{x}_{N_D}^d)\right) &\cdots & \mathcal{L}\left( t_I(\mathbf{x}_{N_D}^d)\right)
 \end{array}
\right],
\
\mathbf{A}_2= \left[
 \begin{array}{ccc}
  \mathcal{B}\left( t_1(\mathbf{x}_1^b)\right) &\cdots & \mathcal{B}\left( t_I(\mathbf{x}_1^b)\right)\\
  \vdots & \ddots & \vdots\\
   \mathcal{B}\left( t_1(\mathbf{x}_{N_{B}}^b)\right) &\cdots & \mathcal{B}\left( t_I(\mathbf{x}_{N_{B}}^b)\right)
 \end{array}
\right],
\end{equation*}
\begin{equation*}
 \mathbf{b} =  \left[
 \begin{array}{c}
  \mathbf{b}_1 \\
  \mathbf{b}_2 \\ 
 \end{array}
\right], \ 
\mathbf{b}_1 = \left[
 \begin{array}{c}
  f(\mathbf{x}_1^d)\\
  \vdots \\
  f(\mathbf{x}_{N_D}^d) \\
 \end{array}
\right],\
\mathbf{b}_2 = \left[
 \begin{array}{c}
  h(\mathbf{x}_{1}^b)\\
  \vdots \\
  h(\mathbf{x}_{N_{B}}^b) \\
 \end{array}
\right],
\end{equation*}
and $\vec{\rho} = \texttt{diag}\left[\sqrt{\frac{\rho_1^d}{N_D}}, \cdots, \sqrt{\frac{\rho_{N_D}^d}{N_D}},\sqrt{\frac{\rho_1^b}{N_{B}}}, \cdots, \sqrt{\frac{\rho_{N_{B}}^b}{N_{B}}}\right]$. Here, $\texttt{diag}$ is a function that creates a diagonal matrix with the elements of the input 1D vector on the diagonal.

Generally, the length of vector $\vec{\alpha}$ is far smaller than the number of collocation points, i.e., $I\ll  N_D+N_B$. Hence, problem \eqref{eq: lstq system} leads to an overdetermined linear system  $\vec{\rho}\mathbf{A} \vec{\alpha} = \vec{\rho} \mathbf{b}$, which can be solved using the least-squares method. Rather than defaulting the penalty weights to 1, we rescale the $L^2$-norm of each row in $\mathbf{A}$ to be equal, i.e., 
\begin{equation}
 \begin{aligned}
  \rho_j^d = \frac{N_D}{ \sum_{i=1}^I \mathcal{L}(t_i( \mathbf{x}_j^d))^2},\ 1\leq j\leq N_D, \ \
  \rho_j^b = \frac{N_{B}}{ \sum_{i=1}^I \mathcal{B}(t_i( \mathbf{x}_j^b))^2} ,\ 1\leq j\leq N_{B}, 
  \label{eq: penalty parameters}
 \end{aligned}
\end{equation}
to balance the contributions of the PDE term and the boundary condition in the loss function. Finally, by incorporating $\vec{\alpha}$ into Eq. \eqref{eq: general form of fto-pinn approximate solution}, the numerical solution of the target linear PDE can be obtained.

\subsubsection{Nonlinear PDEs\label{subsection: nonlinear pde}}
For PDE problem \eqref{eq: general specific pde} with nonlinear operators $\mathcal{L}$ and/or $\mathcal{B}$, substituting $u_{\vec{\alpha}}$ into loss function \eqref{eq: loss function of fto-pinn} yields a nonlinear least-squares problem with respect to $\vec{\alpha}$, given by 
\begin{equation*}
\text{Loss}(\vec{\alpha})= \left\|
  \vec{\rho}\vec{A}(\vec{\alpha})-\vec{\rho}\mathbf{b}
\right\|_2^2,
\end{equation*}
where $\mathbf{b}$ is the same as in Eq. \eqref{eq: lstq system}, while $\vec{A}(\vec{\alpha}):= \left[\vec{A}_1(\vec{\alpha})^T, \vec{A}_2(\vec{\alpha})^T \right]^T$, $\vec{A}_1(\vec{\alpha}) := \left[\mathcal{L}(u_{\vec{\alpha}}(\mathbf{x}_1^d)) , \cdots, \mathcal{L}(u_{\vec{\alpha}}(\mathbf{x}_{N_D}^d)) \right]^T$ and $\vec{A}_2(\vec{\alpha}) := \left[\mathcal{B}(u_{\vec{\alpha}}(\mathbf{x}_1^b)) , \cdots, \mathcal{B}(u_{\vec{\alpha}}(\mathbf{x}_{N_B}^b)) \right]^T$.

To efficiently solve this system, we simply employ the Newton-Linear Least-Squares (Newton-LLSQ) method \cite{dong2021local}. By combining Newton's method and linear least squares, this approach is well-suited for tackling nonlinear problems.
Specifically, we set the output of the branch net to the initial value of $\vec{\alpha}$ and update it iteratively as follows:
\begin{equation*}
    \vec{\alpha}^{(k+1)} = \vec{\alpha}^{(k)} +\delta \vec{\alpha},
\end{equation*}
where the increment $\delta \vec{\alpha}$ is obtained by solving the following linear least-squares problem using the method mentioned in Section \ref{subsection: linear pde}:
\begin{equation*}
        \vec{\rho}\vec{J}(\vec{\alpha}^{(k)})\delta \vec{\alpha} = \vec{\rho}\mathbf{b}-\vec{\rho}\vec{A}(\vec{\alpha}^{(k)}).
\end{equation*}
Here, $\vec{J}(\vec{\alpha}^{(k)})$ is the Jacobian matrix of $\vec{A}(\vec{\alpha})$ evaluated at $\vec{\alpha}^{(k)}$, and $\vec{\rho}$ is determined by rescaling the $L^2$-norm of each row in $\vec{J}(\vec{\alpha}^{(k)})$ to be equal, similar to Eq. \eqref{eq: penalty parameters}. Note that other well-established nonlinear solvers are also applicable.

\subsubsection{Interface problems\label{sec: interface problem}}
Consider an elliptic interface problem defined on a domain $\Omega$ that encloses a subdomain $\Omega_1$:
\begin{equation}
\label{eq: interface problem with complex interface}
\begin{aligned}
\mathcal{L}_i(u(\mathbf{x})):= -\nabla\cdot(a_i\nabla u(\mathbf{x})) &= f_i,\quad \mathbf{x} \in\Omega_i, \ i=1,2,\\
 u(\mathbf{x})|_{\Omega_2}- u(\mathbf{x})|_{\Omega_1} &=g_d,\quad  \mathbf{x}\in \Gamma, \\
a_2\nabla u(\mathbf{x})|_{\Omega_2}\cdot \mathbf{n} - a_1 \nabla u(\mathbf{x})|_{\Omega_1}\cdot \mathbf{n}&=g_n, \quad \mathbf{x}\in \Gamma,\\
u(\mathbf{x})&=h,\quad \mathbf{x} \in \partial\Omega_2, 
\end{aligned}
\end{equation} 
where $\Omega_2:= \Omega\setminus\Omega_1$, interface  $\Gamma:=\partial\Omega_1$.
Generally, the coefficient $a(\mathbf{x})$ in Eq. \eqref{eq: interface problem with complex interface} varies on the two sides of the interface $\Gamma$, thus the solution $u(\mathbf{x})$ and its gradients could be discontinuous across the interface.

Following our previous work \cite{ wu2024solving}, we adapt a DeepONet-type pre-trained IONet to solve the interface problem. Specifically, FTO-PINN fine-tunes the pre-trained IONet locally, and the approximate solution is given by
\begin{equation}
u_{\vec{\alpha}_1, \vec{\alpha}_2}(\mathbf{x}) = \left\{\begin{aligned}
& u_{\vec{\alpha}_1} :=\sum_{i=1}^{I}\alpha_i^1t_i^1(\mathbf{x}),\quad \mathbf{x} \in \Omega_1,\\
& u_{\vec{\alpha}_2} :=\sum_{i=1}^{I}\alpha_i^2t_i^2(\mathbf{x}),\quad \mathbf{x} \in \Omega_2,\\
\end{aligned}\right.
\label{eq: solution form of interface problem}
\end{equation}
where $\vec{\alpha}_j=[\alpha_1^j,\cdots, \alpha_{I}^j]^T\in \mathbb{R}^I$ with $j=1$ and $2$ are new added trainable parameters, and $\{t_i^j\}_{i=1}^{I}$ represents the outputs of the trunk net in IONet restricted to $\Omega_j$. Analogously, by incorporating $u_{\vec{\alpha}_1, \vec{\alpha}_2}(\mathbf{x})$ into the interface problem \eqref{eq: interface problem with complex interface}, the optimization problem \eqref{eq: loss function of fto-pinn} can be rearranged into a linear least-squares problem similar to \eqref{eq: lstq system}. For detailed definitions of $\mathbf{A}$, $\mathbf{b}$ and $\vec{\rho}$, please refer to \ref{appendix:1}.

\subsubsection{Expansion of trunk nets \label{sec: generalization of fto-pinn}}

Inspired by ensemble learning techniques \cite{polikar2012ensemble}, we propose aggregating two or more pre-trained DeepONets to further enhance the performance of FTO-PINN. Given $J\geq 2$ pre-trained DeepONets, and denoting the outputs of the trunk net in the $j$-th DeepONet as $\{t_{1,j},\cdots, t_{I_j,j}\}$, we freeze all pre-trained neural models and fine-tune a combination of the trunk nets by adding a small number of trainable parameters. Then, the approximate solution in FTO-PINN formulated as follows:
\begin{equation}
u_{\vec{\alpha}}(\mathbf{x}) =\sum_{j=1}^{J}\sum_{i=1}^{I_j} \alpha_{ij}t_{i,j}(\mathbf{x}).
\label{eq:generalized FTO-PINN}
\end{equation}
Note that this ensemble technique does not change the linear combination property of trunk nets in (\ref{eq: general form of fto-pinn approximate solution}) and thus the trainable parameter $\vec{\alpha }$ here can be determined as before.

In this work, we restrict our model selection to publicly available pre-trained models. Nevertheless, by carefully extending a single DeepONet, we can still achieve improved performance. Specifically, we propose scaling the input of the trunk net using a set of fixed scale parameters $p_{set}:=\{p_j>0\}_{j=1}^{J}$. The approximate solution in FTO-PINN that incorporates this scaling technique can be expressed as:
\begin{equation}
 u_{\vec{\alpha}}(\mathbf{x})=\sum_{j=1}^{J}\sum_{i=1}^I \alpha_{ij}t_i(p_j\cdot\mathbf{x}),
 \label{eq:generalized FTO-PINN with scales}
\end{equation}
where $p_j\cdot\mathbf{x} := [p_jx_1, \cdots, p_jx_{d_0}]^T$ for any $\mathbf{x}\in\mathbb{R}^{d_0}$.

\subsubsection{Further fine-tune the trunk net as needed.\label{sec: further fine-tuning strategies}}

A prerequisite for the previous algorithm (referred to as phase 1 of FTO-PINN) is that the pre-trained DeepONet has acquired relevant features and patterns for the target PDE. However, this approach may encounter limitations when applied to PDEs with parameter settings that significantly deviate from the distribution of the pre-trained DeepONet's training data.  To mitigate this problem, we introduce phase 2 of FTO-PINN to further fine-tune the trunk net update to ensure it better captures the relevant features of the target PDEs, thereby improving the overall performance of FTO-PINN.  

Trunk net in DeepONet typically employs a FNN \eqref{eq: fully connection network}, which contains multiple linear layers. To adapt to diverse PDEs while maintaining efficiency, we incorporate a low-rank adaptation \cite{lora2022hu} into the FTO-PINN framework to update the weight matrix of the linear layers during the adaptation process. Specifically, for a per-trained linear layer $\mathbf{W}_l\psi_{l-1} + \mathbf{b}_l$ where $\mathbf{W}_l\in \mathbb{R}^{d_l\times d_{l-1}}$ and $\mathbf{b}_l\in \mathbb{R}^{d_l}$, we constrain its update by representing it as:
\begin{equation*}
(\mathbf{W}_l+ \delta \mathbf{W}_l)\psi_{l-1} + \delta\mathbf{b}_l,
\end{equation*}
where $\delta \mathbf{W}_l := BA$ is the product of two low-rank trainable matrices $B\in \mathbb{R}^{{d_l}\times r}$ and $A\in \mathbb{R}^{r\times {d_{l-1}}}$, with $r$ being significantly smaller than $\min(d_l, d_{l-1})$. This modified forward pass results in new outputs for the $l$-th layer: 
\begin{equation*}
 \Tilde{\psi}_l =\sigma( (\mathbf{W}_l+ \delta \mathbf{W}_l)\psi_{l-1} + \delta\mathbf{b}_l).
\end{equation*}
Applying it to selected linear layers, the new outputs for the trunk net are given as
$[\Tilde{t}_1(\mathbf{x}), \Tilde{t}_2(\mathbf{x}), \ldots, \Tilde{t}_I(\mathbf{x})]^T\in\mathbb{R}^I$. 
The approximation function is then formulated as
$u_{\vec{\alpha}} = \sum_{i=1}^I\alpha_i\Tilde{t}_i(\mathbf{x})$. Next, we train the added parameters ${\delta\mathbf{W}_l, \delta\mathbf{b}_l, \vec{\alpha}}$ using the PINN loss (\ref{eq: pinn loss}) for a limited number of iterations. At the beginning of training, $\delta\mathbf{b}_l$ is initialized with the pre-trained values of $\mathbf{b}_l$; $\vec{\alpha}$ is initialized via phase 1 of FTO-PINN; matrices $A$ and $B$ are initialized with random Gaussian values and zeros, respectively, making $\delta \mathbf{W}_l = BA$ initially zero. During training, $\mathbf{W}_l$ remains frozen and does not receive gradient updates.
After training, the updated trunk net captures the relevant features. We then freeze the updated trunk net and optimize $\vec{\alpha}$ using the phase 1 procedures outlined in Sections \ref{subsection: linear pde} and \ref{subsection: nonlinear pde}.

\begin{remark} 
When modifying all linear layers in the pre-trained DeepONet, we can roughly recover the expressiveness of full fine-tuning by setting the rank $r$ to the rank of the pre-trained weight matrices.  
\end{remark}

\subsection{Connection with other numerical methods} 

Whereas the solution space of a PDE with different parameter settings may be infinite-dimensional, there is a long history in computational science of finding a suitable approximate solution of a specific PDE problem by a finite sum
\begin{equation*}
 u(\mathbf{x}) = \sum_{i=1}^I\alpha_i b_i(\mathbf{x}),
\end{equation*}
where $\{b_i\}_{i=1}^I$ are the trial (or basis) functions, and the expansion coefficients $\{\alpha_i\}$ are to be determined. 
The choice of $\{b_i\}_{i=1}^I$ is one of the main features that distinguishes different numerical methods. For instance, $\{b_i\}_{i=1}^I$ in traditional finite element method are functions local in character with finite regularities (e.g., piecewise polynomials), spectral methods employ globally trigonometric functions or orthogonal polynomials as basis functions. Some of meshless methods use local or global radial basis functions, while RWMs adopt shallow NNs with static weights and biases to construct approximate function spaces.

From a function approximation standpoint, FTO-PINN provides a new way to construct the approximation space of the solution. Formally, the output function space of FTO-PINN is set to 
\begin{equation*}
 \mathcal{H} = \text{span}\left\{ t_1(\mathbf{x}), t_2(\mathbf{x}), \cdots, t_I(\mathbf{x}) \right\},
\end{equation*}
and the degrees of freedom (DOF) of $\mathcal{H}$ is $\text{dim}(\mathcal{H})=I$. Then it is straightforward to verify that the following propositions hold.

\begin{proposition}
Suppose that $X$ and $Y$ are Banach spaces, $K \subset X $ are compact set, $\mathcal{G}: K\rightarrow \mathcal{G}(K)\subset Y$ is a continuous operator. For $\varepsilon>0$, assume there exists a DeepONet $\mathcal{G}_{\vec{\theta}}$ of form (\ref{eq: DeepONet}) such that
\begin{equation*}
\sup_{v\in K}\norm{\mathcal{G}(v) - \mathcal{G}_{\vec{\theta}}(v)} \leq \varepsilon.
\end{equation*}
If the function space $\mathcal{F} \subset \mathcal{G}(K)$, then
\begin{equation*}
\sup_{f\in \mathcal{F}}\inf_{\vec{\alpha} \in \R^I}\norm{f(\cdot) - \langle \vec{\alpha}, \ \mathcal{N}_t(\cdot)\rangle } \leq \varepsilon,
\end{equation*}
where $\langle \cdot, \cdot \rangle$ denotes the dot product in $\R^I$. 
\label{pro: one}
\end{proposition}

The DeepONet architecture is motivated by the universal approximation theorem for operators  \cite{chen1995universal}. Proposition \ref{pro: one} shows that an effective pre-trained DeepONet can provide an appropriate set of basis functions to approximate the target solution function. Moreover, as shown in Proposition \ref{pro: two}, FTO-PINN and its generalizations described in Section \ref{sec: generalization of fto-pinn} 
can enhance the accuracy of pre-trained DeepONet models by leveraging a more effective set of basis functions.

\begin{proposition}
Suppose that $X_j, j=1,\cdots, J$ and $Y$ are Banach spaces, $K_j \subset X_j $ with $j=1,\cdots,J$ are compact sets, where $\mathcal{G}_j: K_j\rightarrow \mathcal{G}_j(K_j)\subset Y$ is a continuous operator. For $\varepsilon>0$, assume there exists a DeepONet $\mathcal{G}^j_{\vec{\theta}}$ of form (\ref{eq: DeepONet}) such that
\begin{equation*}
\sup_{v_j\in K_j}\norm{\mathcal{G}_j(v_j) - \mathcal{G}^j_{\vec{\theta}}(v_j)} \leq \varepsilon.
\end{equation*}
Consider the generalized FTO-PINN \eqref{eq:generalized FTO-PINN}, augmented by these $J$ pre-trained DeepONets, defines a linear space denoted by 
\begin{equation}
 \mathcal{H}^e := \text{span}\left\{ t_1^1(\mathbf{x}), \cdots,t_1^{I_1}(\mathbf{x}), \cdots, t_J^1(\mathbf{x}) , \cdots, t_J^{I_J}(\mathbf{x}) \right\},
 \label{eq: he space of fto-pinn}
\end{equation}
with dimension $\text{dim}(\mathcal{H}^e)=\sum_{j=1}^JI_j$,
if the target function space $\mathcal{F} \subset \mathcal{G}(K):=\cup_{j=1}^J\mathcal{G}_j(K_j)$, then we have 
\begin{equation*}
\sup_{f\in \mathcal{F}}\inf_{\substack{\vec{\alpha}_j \in \R^{I_j}\\j=1,\cdots,J}}\norm{f(\cdot) - \sum_{j=1}^J\left\langle \vec{\alpha}_j, \ \mathcal{N}^j_t(\cdot)\right\rangle } \leq \varepsilon.
\end{equation*}
\label{pro: two}
\end{proposition}

\section{Numerical Results\label{sec:numerical results}}
In this section, we perform several experiments to investigate the effectiveness of FTO-PINN, including linear (Section \ref{example 1}), nonlinear (Sections \ref{example 2} and  \ref{example 3}), and interface (Section \ref{example 4}) PDE problems. The baselines for comparison include: (i) PINNs \cite{raissi2019physics, wu2022interfaced}; (ii) DeepONets \cite{lu2021learning, wu2024solving}; and (iii) RWMs \cite{scardapane2017randomness,huang2011extreme, dwivedi2020physics}.

In all test examples, the pre-trained DeepONet-type models used in FTO-PINN are all equipped with smooth activation function \texttt{tanh}. 
Due to the limited publicly available pre-trained DeepONet models, we opted to expand a single pre-trained DeepONet as described in Section \ref{sec: generalization of fto-pinn}. In Eq. \eqref{eq:generalized FTO-PINN with scales}, we select $p_{set}:=\{p_j\}_{j=1}^{J}$ to ensure that coordinate inputs, after the first linear layer, fall within the activation range of the \texttt{Tanh} function (approximately [-3, 3]). This facilitates information propagation through the neural network. 
Specifically, we set $p_{set} := \{1\}$ for $J=1$, $\{1, 0.1\}$ for $J=2$, and 
$p_{set} = \{1\}\cup\{0.1 + jh\}_{j=0}^{J-2}$ for $J>2$,
where $h = \frac{p-0.1}{J-2}$ and $p\in \mathbb{R}^+$ is a constant chosen to approximately satisfy $\|p\cdot \mathbf{W}_1 \mathbf{x}+\mathbf{b}_1\|_\infty\leq \mathbf{3}$ for $\mathbf{x}\in \mathcal{T}_{D}$.

Unless otherwise specified, we employ FTO-PINN without phase 2 in Section \ref{sec: further fine-tuning strategies} by default. The neural network architecture used in both PINN and RWM are FNNs \eqref{eq: fully connection network} equipped with \texttt{Tanh} activation. 
In RWM, the weights and biases of the hidden layer parameters are set to random values uniformly sampled in $[-1,1]$, and the DOF is equal to the width of the FNN output layer. During training, RWM follows the same strategy as FTO-PINN to generate penalty weights \eqref{eq: penalty parameters} and applies the least-squares method to determine the trainable parameters. In contrast, PINN sets all penalty weights to 1 and uses the Adam optimizer with default settings for training. Note that both FTO-PINN and RWM solve the nonlinear examples (Sections \ref{example 2} and   \ref{example 3}) using the one-step Newton-LLSQ iterative method \cite{dong2021local}, as described in Section \ref{subsection: nonlinear pde}. 
To measure the accuracy of neural models, we calculate following errors between the reference solution $u$ and numerical solution $u_{\vec{\alpha}}$:
\begin{itemize}
 \item Maximum absolute point-wise error:
\begin{equation*}
  L^\infty = \max_{1\leq i\leq N}\Big|u(\mathbf{x_i})-u_{\vec{\alpha}}(\mathbf{x_i})\Big|.
\end{equation*}
 \item Relative $L^2$ error:
 \begin{equation*}
  \text{Rel.}\ L^2 = \sqrt{\frac{\sum_{i=1}^N |u(\mathbf{x_i})-u_{\vec{\alpha}}(\mathbf{x_i})|^2}{\sum_{i=1}^N |u(\mathbf{x_i})|^2}}.
 \end{equation*}
\end{itemize}
Here, $N$ is the number of test data points over the entire domain. See \ref{appendix:2} for details on the training point set used for evaluating the loss function \eqref{eq: loss function of fto-pinn} and the test point set used for evaluating the test error.

All experiments are performed on a server equipped with an NVIDIA Tesla V100 GPU and an Intel Xeon E5-2630 v4 CPU. The derivative information for neural models is computed using the \texttt{torch.autograd.grad} function from PyTorch. Overdetermined linear systems, derived from equation \eqref{eq: lstq system}, are solved using the \texttt{linalg.lstsq} function from NumPy \cite{walt2011numpy}.

\subsection{Example 1: Advection equation\label{example 1}}

To demonstrate the effectiveness of FTO-PINN, we begin by evaluating its performance on an advection-dominated PDE explored in \cite{wang2021learning}:
\begin{equation}
  \frac{\partial u(x,t)}{\partial t}+a(x)\frac{\partial u(x,t)}{\partial x} = 0, \;(x,t)\in(0,1)\times(0,1),
  \label{eq: advection equation}
\end{equation}
where the initial and boundary conditions are given as $u(0, t ) = \sin(\frac{\pi}{2}t)$ and 
$u(x, 0) = \sin (\pi x)$, respectively.
Specifically, \cite{wang2021learning} learns a solution operator $\mathcal{G}$ that maps the positive coefficient $a(x)$ to the latent solution $u(x,t)$, where $a(x):= v(x )-\min_x v(x) + 1$ and $v(x)$ sampled from a mean-zero Gaussian random field (GRF) with a radial basis function (RBF) kernel
\begin{equation*}
k_{\beta}(x_1, x_2)=\text{exp}\left(-\frac{\|x_1-x_2\|^2}{2\beta^2}\right) 
\end{equation*}
using a length scale $\beta=0.2$. 
Note that $\beta>0$ controls the smoothness of functions in the GRF, with larger $\beta$ generally producing smoother $v(x)$.

In this example, we fine-tune a pre-trained PI-DeepONet $\mathcal{G}_{\vec{\theta}}$ from \cite{wang2021learning} \footnote{\url{https://github.com/PredictiveIntelligenceLab/Physics-informed-DeepONets}} to solve Eq. \eqref{eq: advection equation}. The branch and trunk nets within $\mathcal{G}_{\vec{\theta}}$ are two modified FNNs, each employing a \texttt{tanh} activation function and consisting of 7 layers with 100 units per layer.
As reported in \cite{wang2021learning}, the final test relative $L^2$ error of $\mathcal{G}_{\vec{\theta}}$ is $2.24\%$. To test the model accuracy,  we numerically solve Equation \eqref{eq: advection equation} using the Lax-Wendroff scheme \cite{iserles2009first} on a fine uniform grid of 513 by 513 points for each test coefficient $a(x)$ to obtain highly accurate reference solutions.

\begin{figure}[htbp] 
	\centering
	\scalebox{0.6}{\includegraphics{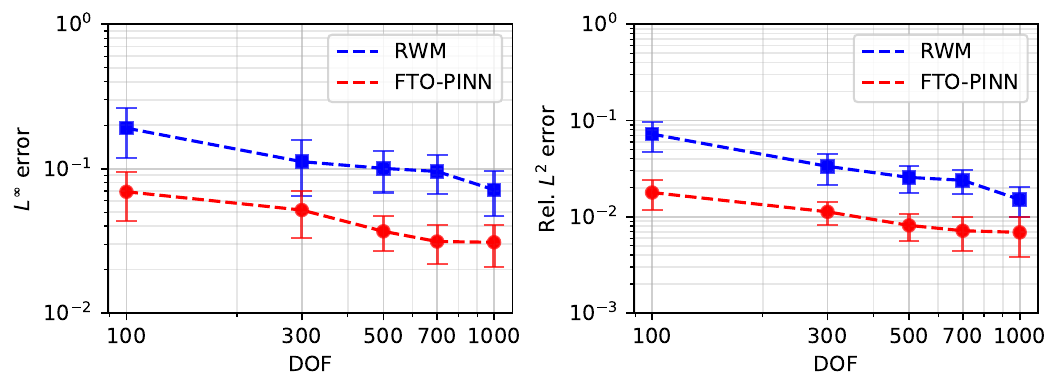}}
	\caption{Advection equation: Variation of the $L^\infty$ and relative $L^2$ errors of FTO-PINN and RWM  with respect to DOF.}
\label{fig: advection equation 2}
\end{figure}

We first discuss the effect of the DOF on the accuracy of FTO-PINN and RWM for solving the advection equation \eqref{eq: advection equation}.  Fig. \ref{fig: advection equation 2} shows the $L^\infty$ and relative $L^2$ errors of these two methods with different DOF for 100 randomly sampled test functions drawn from the GRF with a length scale of $\beta = 0.2$. Here, the results presented for RWM correspond to an FNN with a depth of 2, which achieved the best performance among depths ranging from 1 to 7, with a fixed DOF of 500. As illustrated in the figure, the errors of both FTO-PINN and RWM decrease as DOF increases. Nevertheless, FTO-PINN consistently outperforms RWM, suggesting that fine-tuning a pre-trained DeepONet with PDE-related knowledge within the PINN framework is more effective in capturing the solution of the target equation than employing randomly initialized basis functions.

\begin{table}[htbp]
\centering
\fontsize{8}{7}\selectfont
\begin{threeparttable}
\caption{The relative $L^2$ errors and their comparison for various $a(x)$ generated with different length scales in Advection equation.}
\begin{tabular}{ccccccccccc}
\toprule 
  \multirow{2}{*}{Length scale} & 
 \multicolumn{2}{c}{PI-DeepONet}&
 \multicolumn{2}{c}{PINN}&
 \multicolumn{2}{c}{FTO-PINN (100)$^{\rm a}$} & \multicolumn{2}{c}{FTO-PINN (700)$^{\rm a}$} & \multicolumn{2}{c}{RWM (700)$^{\rm a}$}\cr
 \cmidrule(lr){2-3} \cmidrule(lr){4-5} \cmidrule(lr){6-7} \cmidrule(lr){8-9} \cmidrule(lr){10-11}
 &Rel. $L^2$ & Time(s) & Rel. $L^2$ & Time(s) &Rel. $L^2$ & Time(s) &Rel. $L^2$ & Time(s) &Rel. $L^2$ & Time(s) \cr
 \midrule
$\beta=0.2$ 
& 1.99e-2
& 0.003
& 6.13e-3 & 422
& 1.41e-2 & 0.42
& 5.52e-3 & 1.61
& 2.04e-2 & 2.25
\\
$\beta=0.1$
& 3.03e-2 
& 0.004
& 9.29e-3 & 434
& 2.57e-2 & 0.50 
& 7.60e-3 & 1.62
& 3.36e-2 & 2.23
\\
$\beta=0.05$
& 1.14e-1 
& 0.004
& 5.44e-2 & 430
& 6.85e-2 & 0.47
& 1.54e-2 & 1.63
& 9.20e-2 & 2.29
\\
$\beta=0.025$
& 2.12e-1 
& 0.004
& 3.31e-2 & 419
& 7.13e-2 & 0.49
& 4.24e-2 & 1.64
& 1.98e-1 & 2.23
\\
\bottomrule 
\end{tabular}
\label{tab: advectione-2 equation 1}
\begin{tablenotes}
 \footnotesize
 \item[a] The DOF of the neural models.
 \end{tablenotes}
\end{threeparttable}
\end{table}

\begin{figure}[h] 
	\centering
  \includegraphics[width=1.0\textwidth]{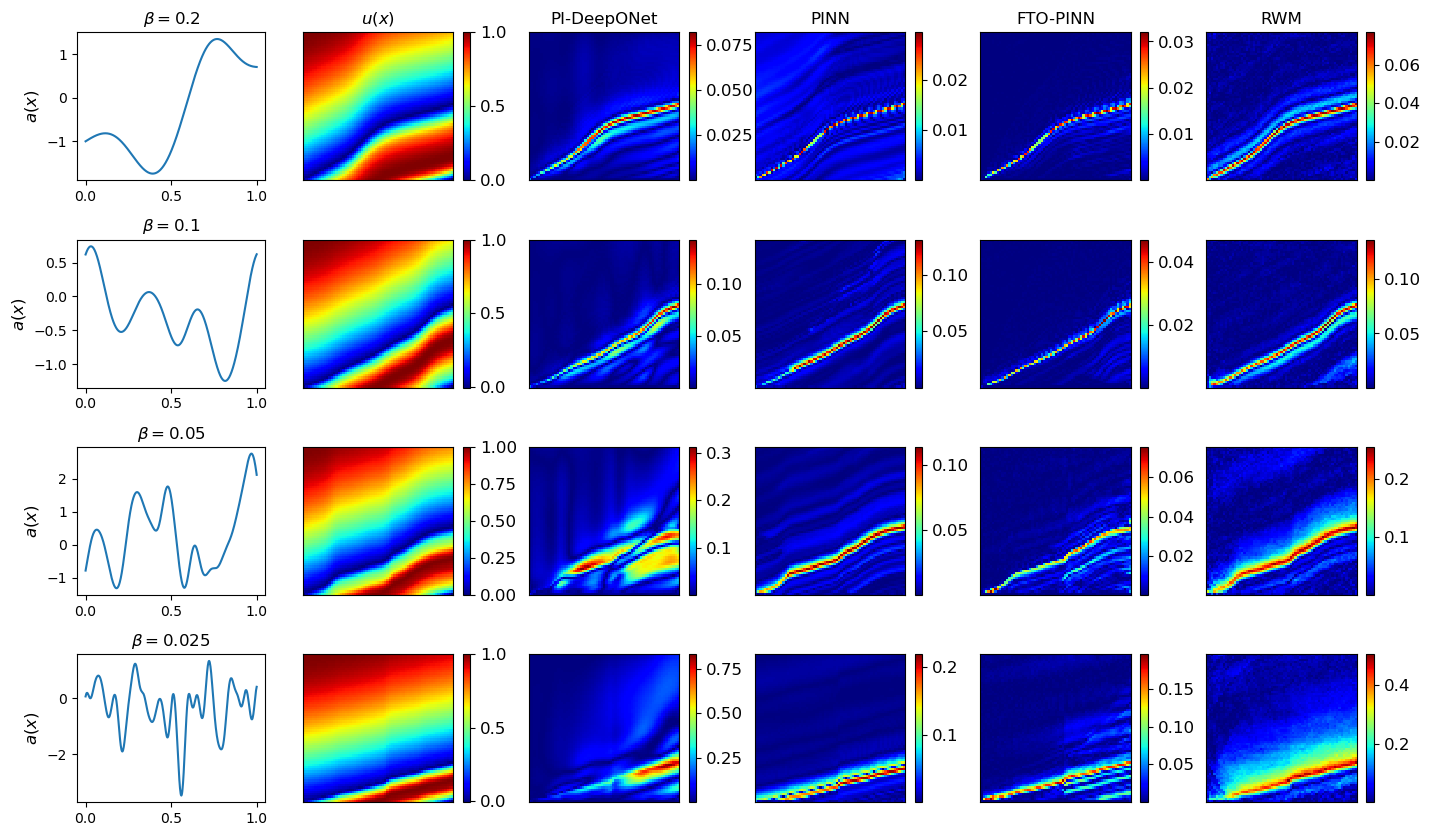}
	\caption{Advection equation: (First column) Coefficient function $a(x)$ generated with different length scales. (Second to sixth columns)  The corresponding reference solutions and point-wise errors for PI-DeepONet, PINN, FTO-PINN (700) and RWM (700).}
\label{fig: advection equation}
\end{figure}

Table \ref{tab: advectione-2 equation 1} compared  relative $L^2$ errors and training costs of pre-trained PI-DeepONet, vanilla PINN, FTO-PINN, and RWM in solving problem \eqref{eq: advection equation} with four different coefficient $a(x)$ (see the first column of Fig. \ref{fig: advection equation} for an illustration). From the setting of hyper-parameters, PINN employs a FNN with 7 layers and 100 units per layer, trained for 50,000 iterations by minimizing the loss function \eqref{eq: pinn loss}. The DOF of RWM is set to 700, while the depth is maintained at 2. As we can see from the table, PI-DeepONet is able to provide fast predictions for specific $a(x)$, but its accuracy is limited, especially for functions $a(x)$ that lie outside the distribution of its training dataset (e.g., those generated with a length scale $\beta < 0.2$). In contrast, PINN exhibits higher accuracy and robustness than PI-DeepONet but requires full model training from scratch, leading to  much longer training times. Without expansion, i,e, with 100 DOF, FTO-PINN demonstrates accuracy that exceeds that of the pre-trained PI-DeepONet. While a slight increase in training time is observed as FTO-PINN increases DOF to 700, this is offset by the significant gains in accuracy, resulting in a favorable balance between accuracy and efficiency compared to PINN. Fig. \ref{fig: advection equation} presents a visual comparison of the point-wise errors obtained by PI-DeepONet, PINN, FTO-PINN, and RWM. Each row depicts a specific coefficient function $a(x)$, with the corresponding reference solution to PDE \eqref{eq: advection equation} in the second column and point-wise error plots for each method in the remaining columns.
Although all methods generally align with the reference, significant errors arise in regions with steep gradients. FTO-PINN demonstrates comparable accuracy to PINN, outperforming DeepONet and RWM in mitigating these errors.

\begin{table}[htbp]
\centering
\fontsize{8}{7}\selectfont
\begin{threeparttable}
\caption{The relative $L^2$ error and training cost of FTO-PINN with further fine-tuned the trunk net. }
\begin{tabular}{cc|cccccccccc}
\toprule 
  \multicolumn{2}{c}{FTO-PINN(700)}& \multirow{2}{*}{\diagbox[]{Rank}{Layer}} & 
 \multicolumn{2}{c}{$l =1^{\rm b}$}&
 \multicolumn{2}{c}{$l =3^{\rm b}$}&
 \multicolumn{2}{c}{$l =5^{\rm b}$} \cr
 \cmidrule(lr){1-2}
 \cmidrule(lr){4-5} \cmidrule(lr){6-7} \cmidrule(lr){8-9} 
 Rel. $L^2$ & Time(s) 
 &
 &Rel. $L^2$ & Time(s) 
 &Rel. $L^2$ & Time(s) &Rel. $L^2$ & Time(s) \cr
 \midrule
\multirow{2}{*}{4.24e-2}&\multirow{2}{*}{1.64}&
$ r = 1$ 
& 3.52e-2 & 29
& 2.86e-2 & 35
& 2.34e-2 & 39
\\
&&
$ r = 5 $
& 3.42e-2 & 31
& 2.68e-2 & 34
& 1.95e-2 & 41
\\
\bottomrule 
\end{tabular}
\label{tab: fine-tune trunk net in example 1}
\begin{tablenotes}
 \footnotesize
 \item[b] Fine-tuning the first $l$ linear layers in the pre-trained trunk-net.
 \end{tablenotes}
\end{threeparttable}
\end{table}

Table \ref{tab: fine-tune trunk net in example 1} records the numerical results of FTO-PINN with phase 2 in Section \ref{sec: further fine-tuning strategies} for solving Eq. \eqref{eq: advection equation}.
In this context, $a(x)$ is generated using $\beta=0.025$, deviating from the distribution of the pre-trained DeepONet's training data, where $\beta=0.2$.
For each case in this table, we fine-tune the modified trunk net by minimizing the PINN loss function $\text{Loss}  = \text{Loss}_{D}  + 10 \text{Loss}_{B}$ after 2000 iterations, with $\text{Loss}_{D}$ and $\text{Loss}_{B}$ representing the loss terms enforcing the PDE and boundary conditions, respectively.  
As shown in the table, FTO-PINN with phase 2 in Section \ref{sec: further fine-tuning strategies} exhibits superior numerical accuracy, even when the rank $r$ is small. Although the fine-tuning process introduces additional computational costs, our method maintains a significant efficiency advantage over the vanilla PINN.

\subsection{Example 2: Diffusion-reaction dynamics\label{example 2}}
To evaluate the effectiveness of FTO-PINN on nonlinear problems, we consider the following diffusion-reaction PDE:
\begin{equation}
 \frac{\partial u}{\partial t} = D\frac{\partial^2 u}{\partial x^2} + k u^2 +f(x), \;(x,t)\in(0,1]\times(0,1]
 \label{eq: nonlinear diffusion-reaction PDE}
\end{equation}
with zero initial and boundary conditions, where the diffusion coefficient is $D=0.01$ and the reaction rate is set to $k=0.01$. This problem configuration is previously examined in \cite{wang2021learning}, where an operator is learned to map the source term $f(x)$ onto the PDE solution $u(x)$. In accordance with the methodology outlined in \cite{wang2021learning}, we approximate the solution operator using a conventional DeepONet model, where the branch and trunk nets are two independent 5-layer \texttt{Tanh} FNNs with 50 units per hidden layer. The training dataset consisted of 10,000 input functions $f(x)$ randomly sampled from a GRF with a length scale of $\beta = 0.2$. For each $f$, the reference solution used in the training process is generated using a second-order implicit finite difference method \cite{iserles2009first} on a 100 by 100 equispaced grid, while that for testing here are generated on a finer grid of 513 by 513 equispaced points. For other specific parameter settings, please refer to the supplementary material in \cite{wang2021learning}.

\begin{figure}[ht] 
	\centering
	\scalebox{0.45 }{\includegraphics{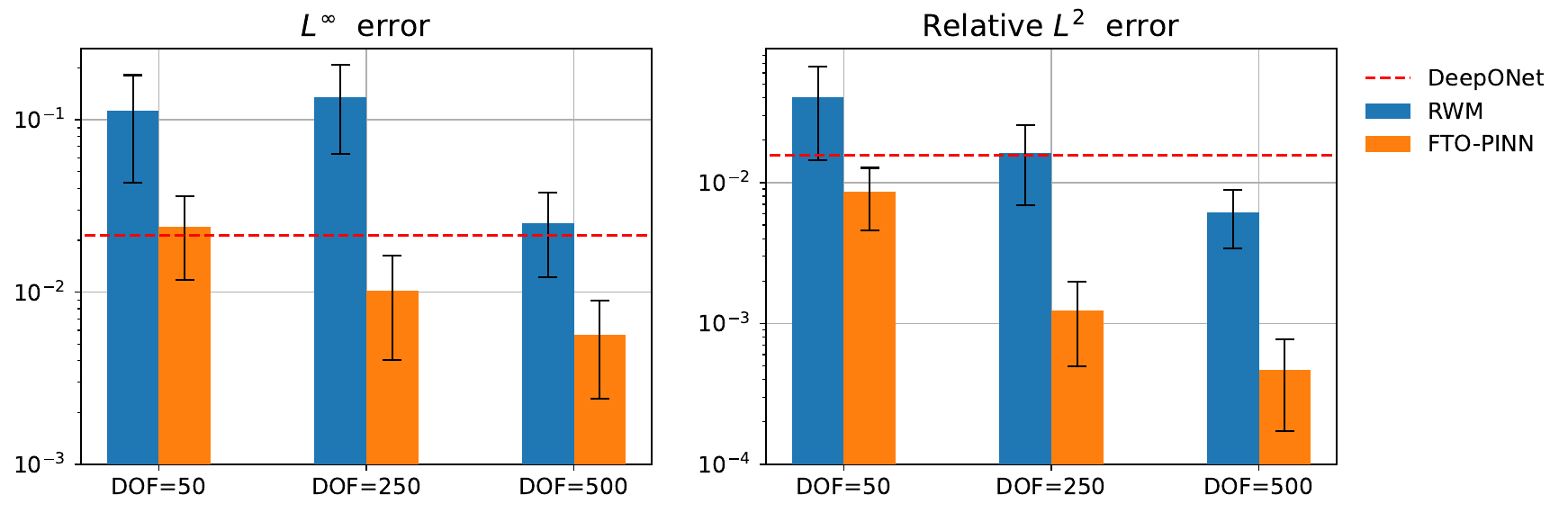}}
	\caption{Diffusion-reaction equation: The $L^\infty$ and relative $L^2$ errors of FT0-PINN and RWM over 100 randomly sampled $f(x)$. Here, red dashed lines indicate the mean error of pre-trained DeepONet.}
\label{fig: diffusion-reaction dynamics 1}
\end{figure}

Fig. \ref{fig: diffusion-reaction dynamics 1} shows the errors of FTO-PINN and RWM for solving the nonlinear PDE \eqref{eq: nonlinear diffusion-reaction PDE} via the one-step Newton-LLSQ iterative method  \cite{dong2021local}. Here, $f(x)$ is generated from a GRF with a length scale of $\beta=0.2$, consistent with the distribution used to create the training set for the pre-trained DeepONet. The results presented for RWM correspond to an FNN with a depth of 3, which achieved the best performance among depths ranging from 1 to 5, with a fixed DOF of 500.
The general trend shows increased prediction accuracy as DOF increases. Building on our success with linear problems, we further show that FTO-PINN outperforms RWM on nonlinear problems, highlighting the effectiveness of FTO-PINN.

\begin{table}[ht]
\centering
\fontsize{8}{7}\selectfont
\begin{threeparttable}
\caption{The relative $L^2$ errors and their comparison for four $f(x)$ generated with different length scales in Diffusion-reaction equation. }
\begin{tabular}{ccccccccc}
\toprule 
  \multirow{2}{*}{Length scale} & 
 \multicolumn{2}{c}{ DeepONet}&
 \multicolumn{2}{c}{PINN}
  & \multicolumn{2}{c}{FTO-PINN (500)} & \multicolumn{2}{c}{RWM (500)}\cr
 \cmidrule(lr){2-3} \cmidrule(lr){4-5} \cmidrule(lr){6-7} \cmidrule(lr){8-9} 
 &Rel. $L^2$ & Time(s) &Rel. $L^2$ & Time(s) &Rel. $L^2$ & Time(s) &Rel. $L^2$ & Time(s) \cr
 \midrule
$\beta =0.2$ 
& 1.29e-2 & 0.001
& 4.39e-2 & 403 
& 4.47e-4 & 1.01
& 6.56e-3 & 3.73
\\
$\beta =0.1$
& 8.51e-2 & 0.001
& 4.93e-2 & 435
& 2.37e-4 & 1.08
& 1.45e-2 & 3.75
\\
$\beta =0.05$
& 1.76e-1 & 0.001
& 5.34e-2 & 489
& 3.54e-3 & 1.07
& 5.00e-2 & 3.70
\\
$\beta =0.025$
& 3.65e-1
& 0.002
& 5.89e-2
& 412
& 3.48e-2 & 1.05
& 1.48e-1 & 3.71 
\\
\bottomrule 
\end{tabular}
\label{tab: diffusion-reaction dynamics 1}
\end{threeparttable}
\end{table}

\begin{figure}[ht] 
	\centering
 \includegraphics[width=1.0\textwidth]{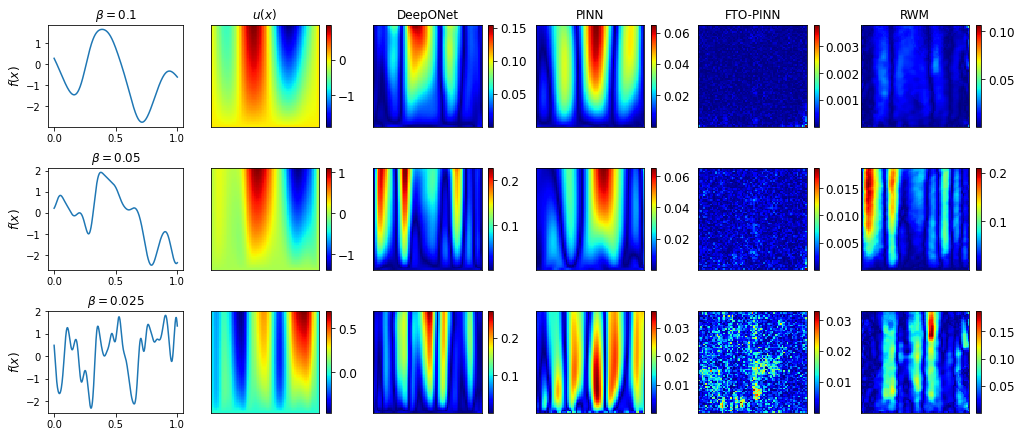}
	\caption{Diffusion-reaction equation: (First column) Source terms $f(x)$ generated with different length scales. (Second to sixth columns) The reference solutions $u(\mathbf{x})$ and point-wise errors for PI-DeepONet, PINN, FTO-PINN (500) and RWM (500).}
\label{fig: diffusion-reaction dynamics 2}
\end{figure}

Table \ref{tab: diffusion-reaction dynamics 1} summarizes the accuracy and training costs of pre-trained DeepONet, PINN, FTO-PINN, and RWM for solving problem \eqref{eq: nonlinear diffusion-reaction PDE}. Here, PINN employs a 5-layer \texttt{Tanh} FNN with 50 units per hidden layer, whereas RWM maintains a depth of 3. The overall trend shows that the accuracy of neural models declines as the length scale $\beta$ decreases. Specifically, when the $\beta$ decreases from 0.2 to 0.025, the accuracy of DeepONet declines by an order of magnitude, dropping from $10^{-2}$ to approximately $10^{-1}$. 
In contrast, PINN, trained for 50,000 iterations to minimize the loss function \eqref{eq: pinn loss}, consistently maintains relative $L^2$ errors around $10^{-2}$ across different length scales. This robustness comes with the trade-off of longer training times due to the optimization of a highly non-convex loss function for solving PDEs. FTO-PINN and RWM improve the accuracy of DeepONet through a one-step nonlinear iteration. 
In all cases, FTO-PINN provides more accurate solutions than DeepONets and RWM, while maintaining comparable accuracy to PINN and significantly reducing computational time. Furthermore, Fig. \ref{fig: diffusion-reaction dynamics 2} compares the reference and predicted solutions for $\beta<0.2$. As shown in the first column, decreasing the value of $\beta$ results in an increase in the number of inflection points produced by $f$. Despite the complexity of the source term $f$, the FTO-PINN predictions consistently show superior alignment with the reference solutions. These findings highlight the capability of FTO-PINN to enhance the efficiency of PINN and improve the accuracy of DeepONet in handling $f$ beyond the distribution of its training data, even in nonlinear scenarios.

\subsection{Example 3: Burgers’ equation \label{example 3}}
In this section, we demonstrate that FTO-PINN can adapt a pre-trained DeepONet, trained to approximate the solution operator of related PDEs, to the target PDE. Consider Burgers' equation with different values of viscosity $\mu$:
\begin{equation}
\begin{aligned}
 \frac{du}{dt}+u\frac{du}{dx}-\mu\frac{d^2u}{dx^2} &= 0, \ (x,t)\in (0,1)\times(0,1],\\
 u(x,0) &= u_0(x), \ x\in (0,1),
 \label{eq:Burgers' equation}
\end{aligned}
\end{equation}
subject to periodic boundary conditions
\begin{align*}
 u(0,t)&=u(1,t),\\
 \frac{du}{dx}(0,t) &= \frac{du}{dx}(1,t).
\end{align*}
For viscosity $\mu=0.01$, a PI-DeepONet model $\mathcal{G}_{\vec{\theta}}$\footnote{\url{https://github.com/PredictiveIntelligenceLab/Physics-informed-DeepONets}} is pre-trained in \cite{wang2021learning} to learn the solution operator that maps the initial condition $u_0(x)$ to the solution $u(x, t)$, where the input function $u_0(x)$ is generated from a GRF $\mathcal{N}(0, 25^2(-\Delta+5^2I)^{-4})$ satisfying the periodic boundary conditions and the trunk net employs a modified FNN with a depth of 7 and 100 units per layer. For each test function in this work, the reference solution is generated on a grid of 201 by 201 equispaced points, while other specific parameter settings please refer to the supplementary material of \cite{wang2021learning}.

We first investigated the performance of pre-trained PI-DeepONet, PINN, FTO-PINN, and RWM in solving Eq. \eqref{eq:Burgers' equation} with $\mu=0.01$. Here, three test initial conditions are generated from different GRFs. In each case, PINN employs a 7-layer \texttt{tanh} FNN with 100 units per layer, trained for 50,000 iterations. The FNN used in RWM has a depth of 2 (selected from a range of 1 to 7) and a width of 500. In addition to using automatic differentiation \cite{baydin2017automatic} to calculate the derivatives of the loss function \eqref{eq: loss function of fto-pinn}, we also explore the use of the central difference method \cite[Section 5]{shang2022deep} with a spatiotemporal step size of $5\times10^{-4}$. The resulting model is referred to as FTO-PINN-FD.

\begin{table}[htbp]
\centering
\fontsize{8}{7}\selectfont
\begin{threeparttable}
\caption{The relative $L^2$ errors and their comparison for three initial conditions $u_0(x)$ generated from different GRFs in Burgers’ equation.
}
\begin{tabular}{lcccccccc}
\toprule 
  \multirow{2}{*}{\diagbox[]{Method}{GRF}} & 
 \multicolumn{2}{c}{$\mathcal{N}(0, 25^2(-\Delta+5^2I)^{-3.5})$ }&
 \multicolumn{2}{c}{$\mathcal{N}(0, 25^2(-\Delta+2^2I)^{-4})$}&
  \multicolumn{2}{c}{$\mathcal{N}(0, 5^2(-\Delta+5^2I)^{-4})$} \cr
 \cmidrule(lr){2-3} \cmidrule(lr){4-5} \cmidrule(lr){6-7} 
 &Rel. $L^2$ & Time(s) & Rel. $L^2$ & Time(s) &Rel. $L^2$ & Time(s) \cr
 \midrule
PI-DeepONet 
& 1.72e-1 & 0.002
& 6.43e-3 & 0.002 
& 2.04e-2 & 0.002

\\
PINN
& 9.87e-2 & 877
& 8.54e-2 & 859 
& 6.93e-2 & 832

\\
FTO-PINN (500) 
& 5.22e-2 & 2.18
& 2.72e-4 & 2.05
& 9.59e-5 & 2.02
\\
FTO-PINN-FD (500) 
& 5.86e-2 & 0.49
& 7.82e-3 & 0.41
& 1.25e-3 & 0.44
\\
RWM (500)
& 1.99e-1 & 1.93
& 3.56e-3 & 1.87
& 1.64e-2 & 1.81
\\
\bottomrule 
\end{tabular}
\label{tab: burgers equation 1}
\end{threeparttable}
\end{table}

\begin{figure}[ht] 
	\centering
 \includegraphics[width=1.0\textwidth]{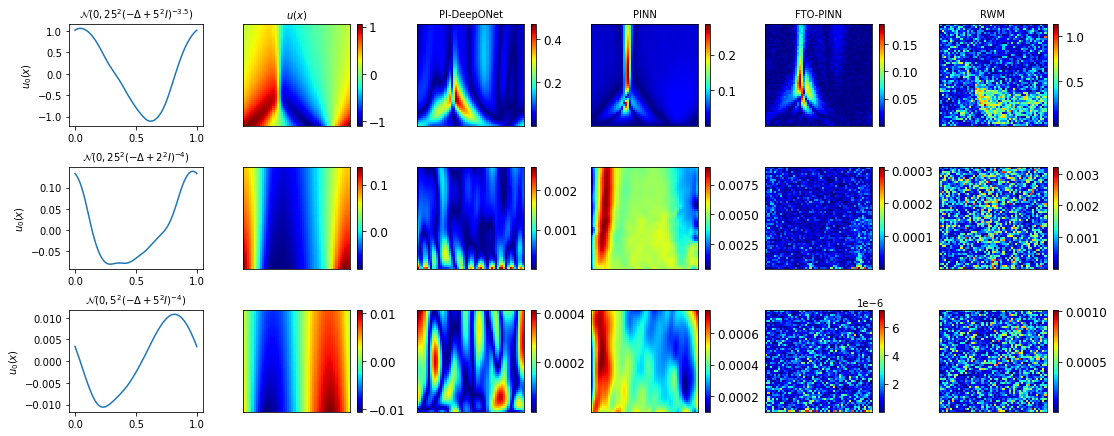}
	\caption{Burgers’ equation: (First column) Initial condition $u_0$ generated with different GRFs. (Second to sixth columns) The reference solutions $u(\mathbf{x})$ and point-wise errors for PI-DeepONet, PINN, FTO-PINN (500) and RWM (500).}
\label{fig:Burgers' equation1}
\end{figure}

Table \ref{tab: burgers equation 1} presents the accuracy and training costs of these neural models for solving Burgers’ equation with $\mu=0.01$. Here, initial conditions $u_0$ are sampled from GRFs that are different from the distribution of the training data for the pre-trained DeepONet $\mathcal{G}_{\vec{\theta}}$. 
As shown in this table, FTO-PINN significantly enhances the generalization performance of the pre-trained DeepONet $\mathcal{G}_{\vec{\theta}}$.  
Compared to PINN, FTO-PINN achieves superior accuracy while demanding substantially fewer computational resources. In addition, FTO-PINN-FD offers a trade-off between training speed and accuracy. It achieves accelerated training times but may introduce errors due to the derivative estimation within the loss function. 
Visualizations of the errors are provided in Fig. \ref{fig:Burgers' equation1}. In all cases, FTO-PINN achieves the lowest point-wise error. However, discrepancies may arise in cases where the solution exhibits steep gradients, shown in the first row. As illustrated in Example 1, specialized approximation strategies or loss functions are necessary to address this challenge. We intend to explore this topic further in future research.

\begin{table}[ht]
\centering
\fontsize{8}{7}\selectfont
\begin{threeparttable}
\caption{The relative $L^2$ errors and their comparison for different values of viscosity in Burgers’ equation.}
\begin{tabular}{cccc}
\toprule 
 Viscosity  $\mu$ & PI-DeepONet& FTO-PINN (500) & RWM (500) \cr
 \midrule
$0.02$         
& 2.17e-1 $\pm$ 1.94e-2
& 7.77e-3 $\pm$ 4.33e-3
& 6.83e-2 $\pm$ 9.52e-2 
\\
$0.01$ 
& 1.08e-2 $\pm$ 1.66e-2
& 9.12e-4 $\pm$ 3.35e-3
& 1.07e-2 $\pm$ 1.58e-2 
\\
$0.005$ 
& 1.34e-1 $\pm$ 2.69e-2
&1.37e-2 $\pm$ 1.89e-2
&1.08e-1 $\pm$ 9.12e-2
 
\\
$0.0025$ 
&2.06e-1 $\pm$ 2.49e-2
&4.60e-2 $\pm$ 4.17e-2
&1.79e-1 $\pm$ 1.19e-1 
\\
\bottomrule 
\end{tabular}
\label{tab:Burgers' equation 2}
\end{threeparttable}
\end{table}

Next, we consider the Burgers’ equation \eqref{eq:Burgers' equation} with varying viscosity values. For each viscosity value, we randomly sampled 100 initial conditions $u_0(x)$ from the GRF $\mathcal{N}(0, 25^2(-\Delta+5^2I)^{-4})$. Table \ref{tab:Burgers' equation 2} reports the $L^2$ errors and compares the results for the pre-trained PI-DeepONet, FTO-PINN, and RWM. Recall that the pre-trained PI-DeepONet is trained for the Burgers’ equation with $\mu=0.01$. As expected, PI-DeepONet achieves its lowest relative $L^2$ test error for $\mu=0.01$, while FTO-PINN achieves at least a tenfold reduction in this error. Furthermore, for other values of viscosity $\mu$, the accuracy of PI-DeepONet significantly decreases to the order of $10^{-1}$. In contrast, FTO-PINN and RWM, utilizing a one-step Newton-LLSQ iterative method, consistently achieve higher prediction accuracy, especially for FTO-PINN. These results demonstrate the effectiveness of FTO-PINN utilizing only PDE-related pre-trained DeepONet models.

\subsection{Example 4: Elliptic interface problem\label{example 4}}
\begin{figure}[htbp] 
	\centering
	\scalebox{0.4}{\includegraphics{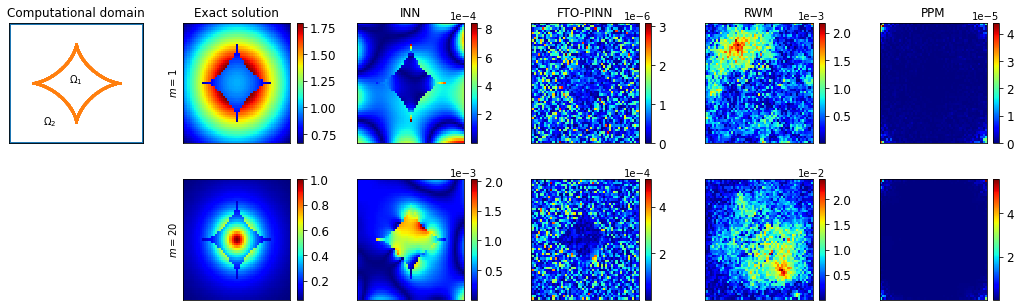}}
	\caption{Elliptic interface problem: (First column)  Computational domain. (Second to sixth columns)  Exact solution and point-wise errors of INN, FTO-PINN (1000), RWM (1000) and PPM (1056).} 
 \label{fig: point-wise error of interface peoblem with complex}
\end{figure}

The proposed FTO-PINN is essentially a meshless method, similar to PINN. In this case, we further investigate the performance of FTO-PINN for solving the elliptic interface problem \eqref{eq: interface problem with complex interface} with an irregular interface $\Gamma$ (see the left panel of Fig. \ref{fig: point-wise error of interface peoblem with complex} for an illustration), which is given as 
\begin{equation*}
 (x_1,x_2)=(0.65\text{cos}(\vartheta)^3,0.65\text{sin}(\vartheta)^3), \ \vartheta\in [0,2\pi).
\end{equation*} 
The computational domain of Eq. \eqref{eq: interface problem with complex interface} is defined as $\Omega=[-1,1]^2$. Here, $\Omega_1$ is a subset of $\Omega$ bounded by the interface $\Gamma$, and $\Omega_2=\Omega/(\Omega_1\cup\Gamma)$. The discontinuous coefficient $a$ is set to
\begin{equation*}
a(x_1,x_2)=\left\{
\begin{aligned}
&2, \ &\text{in} \ \Omega_1, \\
&1, \ &\text{in} \ \Omega_2. \\
\end{aligned}
\right.
\end{equation*}
As solutions to interface problems usually have low global regularity, we consider an exact solution of Eq. \eqref{eq: interface problem with complex interface}, which is discontinuous across the interface, in the following form:
\begin{equation}
u(x_1,x_2)=\left\{
\begin{aligned}
&\frac{1}{1+m(x_1^2+x_2^2)}, \ &\text{in} \ \Omega_1, \\
&\frac{2}{1+m(x_1^2+x_2^2)}, \ &\text{in} \ \Omega_2, \\
\end{aligned}
\right.
\label{eq: solution of interface problem with complex interface}
\end{equation}
where $m>0$ is a given constant. Note that
\begin{equation*}
\frac{\partial^{2k} u}{\partial x_i^{2k}}\Big|_{(0,0)}=(-m)^k(2k)!.
\end{equation*}
Thus, a large $m$ corresponds to a solution with high derivative at $(0,0)$. For any given $m$, the corresponding source terms, boundary conditions, and interface conditions are different and can be derived from Eq. \eqref{eq: solution of interface problem with complex interface}. This problem is also considered in \cite{wu2024solving}, where a pre-trained IONet $\mathcal{G}_{\vec{\theta}}$ \eqref{eq: structure of ionet} is learned to map the source term $f(\mathbf{x})$ to the solution  $u(\mathbf{x})$. The branch and trunk nets in $\mathcal{G}_{\vec{\theta}}$ are \texttt{tanh} FNNs with a depth of 5 and a width of 50. Notably,  Eq. \eqref{eq: solution of interface problem with complex interface} with $m=10$ is included in the solution space of the training data for the pre-trained IONet $\mathcal{G}_{\vec{\theta}}$.

The approximate solution of FTO-PINN, utilizing the pre-trained IONet $\mathcal{G}_{\vec{\theta}}$, takes the form given in Eq. \eqref{eq: solution form of interface problem}. To ensure a fair comparison, RWM employs two random-weighted \texttt{tanh} FNNs, each responsible for solving in a separate subdomain. Similarly, we consider using Pascal polynomials for local approximation of the solution to the interface problem, referred to as PPM \cite{orucc2021efficient}, while other operations remain consistent with RWM. Furthermore, we investigate the performance of INN  \cite{wu2022interfaced}, a PINN-type method that employs an FNN with 5 layers and 50 units per layer, trained over 25,000 iterations by minimizing the loss function \eqref{eq: pinn loss}  derived from Eq. \eqref{eq: interface problem with complex interface}.

\begin{figure}[htb] 
	\centering
	\scalebox{0.43}{\includegraphics{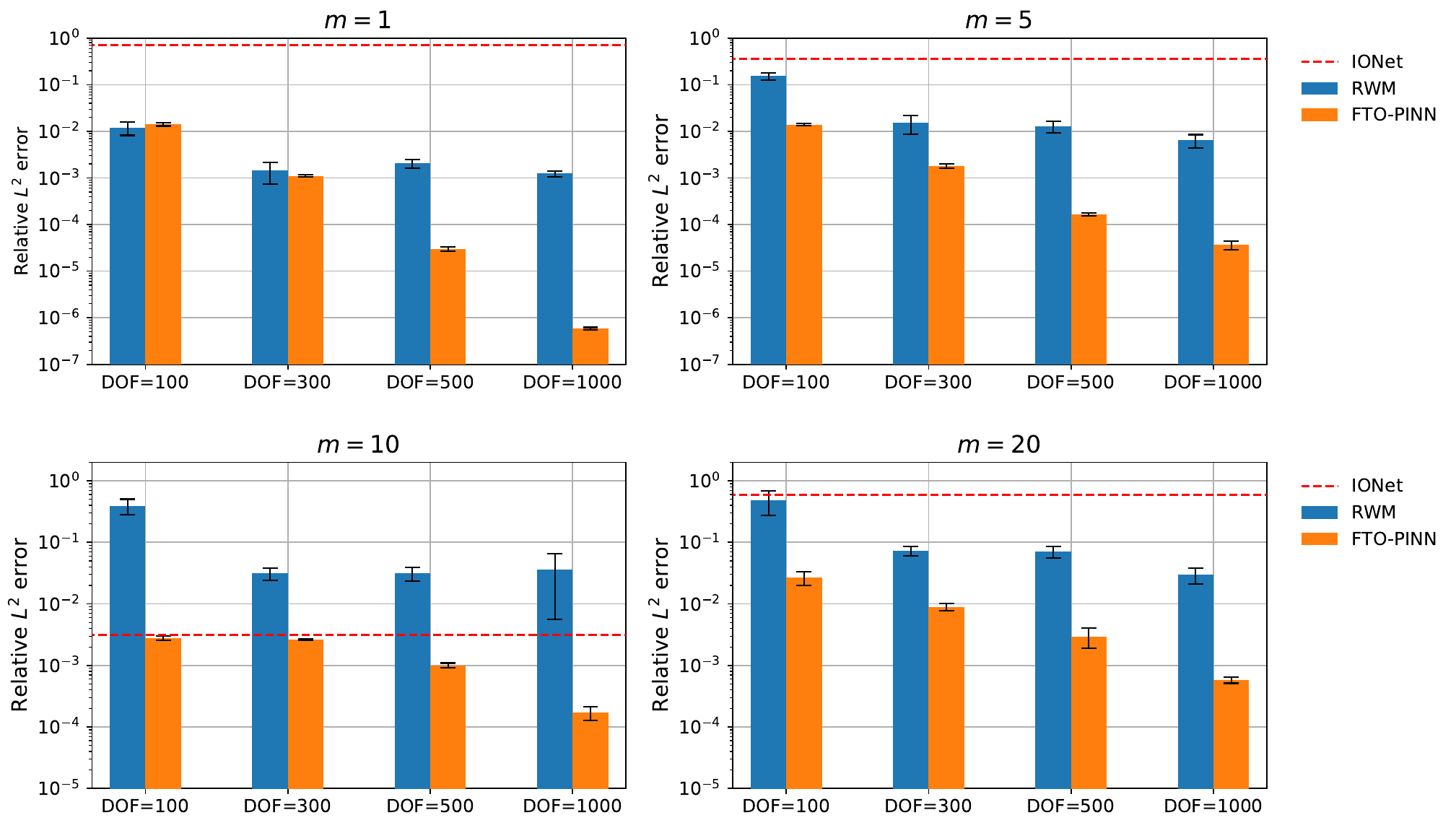}}
	\caption{Elliptic interface problem: Variation of the relative $L^2$ errors of FTO-PINN and RWM with respect to DOF. Here, the dashed lines denote the mean error of IONet.} 
 \label{fig: error bar of interface problem with irregular interfac}
\end{figure}

We initially investigate the influence of DOF on the numerical accuracy of RWM and FTO-PINN for solving the interface problem. Fig. \ref{fig: error bar of interface problem with irregular interfac} illustrates the mean and standard deviation of the relative $L^2$ error across five independent runs. Each run uses a random initialization of the FNN and a different set of points for evaluate the loss function. In this analysis, the depth of RWM is set to 2, chosen from a range of 1 to 5. As illustrated in the figure, the relative $L^2$ error of FTO-PINN generally decreases with increasing DOF. Compared to IONet and RWM, FTO-PINN exhibits superior accuracy and is less sensitive to the value of $m$.

\begin{table}[htbp]
\centering
\fontsize{8}{7}\selectfont
\begin{threeparttable}
\caption{The relative $L^2$ errors and their comparison for different $m$ in the elliptic interface problem.}
\begin{tabular}{lcccccccccc}
\toprule 
  \multirow{2}{*}{\diagbox[]{Method}{$m$}} & 
 \multicolumn{2}{c}{$m=1$ }&
 \multicolumn{2}{c}{$m=5$}&
  \multicolumn{2}{c}{$m=10$}&
  \multicolumn{2}{c}{$m=20$} \cr
 \cmidrule(lr){2-3} \cmidrule(lr){4-5} \cmidrule(lr){6-7} \cmidrule(lr){8-9} 
 &Rel. $L^2$ & Time(s) & Rel. $L^2$ & Time(s) &Rel. $L^2$ & Time(s) &Rel. $L^2$ & Time(s) \cr
 \midrule
IONet 
& 7.19e-1 & 0.002
& 3.57e-1 & 0.002
& 3.16e-3 & 0.002
& 5.90e-1 & 0.002
\\
INN
& 2.67e-4 & 1141
& 5.42e-4 & 1160
& 1.76e-3 & 1128
& 9.38e-3 & 1154
\\
FTO-PINN (1000) 
& 5.68e-7 & 0.87
& 3.66e-5 & 0.75
& 1.50e-4 & 0.74
& 5.94e-4 & 0.76
\\
RWM (1000)
& 5.03e-4 & 3.43
& 7.40e-3 & 3.33
& 7.73e-3 & 3.37
& 3.04e-2 & 3.37
\\
PPM (1056)$^{\rm c}$
& 1.94e-6 & 1.18
& 1.77e-2 & 1.17
& 1.19e-1 & 1.19
& 5.37e-1 & 1.17
\\
\bottomrule 
\end{tabular}
\label{tab: interface problem}
\begin{tablenotes}
 \footnotesize
 \item[c] Training costs of PPM are recorded as CPU time.
 \end{tablenotes}
\end{threeparttable}
\end{table}

Table \ref{tab: interface problem} records relative $L^2$ error and training cost of the pre-trained IONet, INN, FTO-PINN, RWM, and PPM for solving interface problems with different $m$. As shown in the table, for $m = 1$, the relative  $L^2$ error of FTO-PINN reaches the order of $10^{-7}$, outperforming other methods. As the value of $m$ increases, the accuracy of all neural methods generally decreases due to the intensifying singularity of the solution \eqref{eq: solution of interface problem with complex interface} near the origin.  Despite this inherent challenge, FTO-PINN consistently demonstrates its superiority, achieving the lowest relative $L^2$ error with a comparatively minimal investment in training time across all scenarios investigated. Moreover, Fig. \ref{fig: point-wise error of interface peoblem with complex} presents a comparison of point-wise errors of the numerical solutions obtained by the INN, FTO-PINN, RWM, and PPM. As shown in this figure, the behavior of the point-wise error aligns with that of the relative $L^2$ error. FTO-PINN demonstrates better robustness in handling different values of $m$. These observations  highlight the advantages of FTO-PINN in addressing irregular interface problems.

\section{Conclusion\label{sec:conclusion}}
In this work, we propose FTO-PINN, a novel hybrid approach that combines the advantages of PINNs and pre-trained DeepONets for efficient, meshless PDE solving. By leveraging a pre-trained DeepONet that has already learned many relevant features, we can bypass the initial training stages of the PINN and focus on adapting the pre-trained DeepONet model to the specific PDE at hand. 
This significantly reduces the computational time and resources required, improving efficiency compared to standard PINNs. Additionally, FTO-PINN can be regarded as a post-processing technique for DeepONets that are trained on extensive function data to learn operators. Further adaptation to specific PDEs allows the model to retain learned knowledge while leads to improved performance on our specific task. To enhance the adaptability of FTO-PINN, we also explore trunk network expansion and low-rank adaptation techniques. These methods enable the model to effectively handle new PDEs, including those outside the domain for which the DeepONet was initially trained and performs adequately. Numerical experiments on a variety of PDE problems demonstrate the superior performance of FTO-PINN compared to existing baselines.

Herein, due to the requirement of pre-trained DeepONet, we focus primarily on the meshless advantages of PINNs and demonstrate the effectiveness of FTO-PINN across relatively low-dimensional PDEs. The performance of FTO-PINN on high-dimensional problems is left for future investigation.
Future extensions of FTO-PINN could incorporate other pre-trained neural operators and alternative loss construction methods. For example, we plan to modify the fine-tuning scheme to accommodate neural operators such as Fourier neural operator \cite{li2021fourier} or Meta-Mgnet \cite{chen2022meta} in future work. Additionally, while we currently enforce the strong form of PDEs, which is straightforward to implement using automatic differentiation, exploring weak or variational forms may offer further benefits.
Finally, in this study, we assume that the computational domain of the target PDE aligns with the given pre-trained DeepONet. In future work, we intend to explore domain decomposition and scaling strategies to leverage multiple pre-trained neural models for efficiently solving PDEs defined on more complex and diverse domains.

\section*{Acknowledgments}

The author would like to thank Professor Jun Hu for valuable discussions.

\bibliographystyle{elsarticle-num}
\bibliography{ref}

\appendix
\section{Definitions of the linear least-squares problem in Section \ref{sec: interface problem}\label{appendix:1}}

Recall that the interface problem has the following form:
\begin{subequations} 
\begin{align}
\mathcal{L}_i(u(\mathbf{x})):= -\nabla\cdot(a_i\nabla u(\mathbf{x})) &= f_i,\quad \mathbf{x} \in\Omega_i, \ i=1,2  \label{eq:1A},\\
 u(\mathbf{x})|_{\Omega_2}- u(\mathbf{x})|_{\Omega_1} &=g_d,\quad  \mathbf{x}\in \Gamma, \label{eq:1B}, \\
a_2\nabla u(\mathbf{x})|_{\Omega_2}\cdot \mathbf{n} - a_1 \nabla u(\mathbf{x})|_{\Omega_1}\cdot \mathbf{n}&=g_n, \quad \mathbf{x}\in \Gamma \label{eq:1C},\\
u(\mathbf{x})&=h,\quad \mathbf{x} \in \partial\Omega_2 \label{eq:1D}.
\end{align}
\end{subequations}
Let $\{\mathbf{x}^{d_k}_j\}_{j=1}^{n_{d_k}}$ with $k=1,2$, $\{\mathbf{x}^{\gamma}_j\}_{j=1}^{n_\gamma}$and $\{\mathbf{x}^{b}_j\}_{j=1}^{n_{b}}$ be the sets of randomly collected points from the domains $\Omega_1$ and $\Omega_2$, the interface $\Gamma$, and the boundary $\partial\Omega_2$, respectively.  By incorporating  $u_{\vec{\alpha}_1,\vec{\alpha}_2}$ \eqref{eq: solution form of interface problem} into the interface problem, the loss function has the following specific form:
\begin{equation*}
\begin{aligned}
\text{Loss}(\vec{\alpha}_1, \vec{\alpha}_2)
&=\frac{1}{n_{d_1}}\sum_{j=1}^{n_{d_1}}\rho_j^{d_1}\left| \mathcal{L}_1\left( u_{\vec{\alpha}_1}(\mathbf{x}_j^{d_1})\right)-f_1(\mathbf{x}_j^{d_1})\right|^2
+ \frac{1}{n_{d_2}}\sum_{j=1}^{n_{d_2}}\rho_j^{d_2}\left| \mathcal{L}_2\left( u_{\vec{\alpha}_2}(\mathbf{x}_j^{d_2})\right)-f_2(\mathbf{x}_j^{d_2})\right|^2\\ &+\frac{1}{n_{b}}\sum_{j=1}^{n_{b}}\rho_j^b \left|  u_{\vec{\alpha}_2}(\mathbf{x}_j^b)-h(\mathbf{x}_j^b)\right|^2 +
\frac{1}{n_{\gamma}}\sum_{j=1}^{n_{\gamma}}\rho_j^{\gamma_d} \left|  u_{\vec{\alpha}_2}(\mathbf{x}_j^\gamma)- u_{\vec{\alpha}_1}(\mathbf{x}_j^\gamma)-g_d(\mathbf{x}_j^\gamma)\right|^2
\\
&+ \frac{1}{n_{\gamma}}\sum_{j=1}^{n_{\gamma}}\rho_j^{\gamma_n} \left|  a_2 \nabla u_{\vec{\alpha}_2}(\mathbf{x}_j^\gamma)\cdot \mathbf{n} - a_1 \nabla u_{\vec{\alpha}_1}(\mathbf{x}_j^\gamma)\cdot \mathbf{n}-g_n(\mathbf{x}_j^\gamma)\right|^2,
\end{aligned}
\end{equation*}
where the first two terms enforce the PDE constraints \eqref{eq:1A}, the third term penalizes violations of the boundary condition \eqref{eq:1D}, and the fourth and fifth terms enforce interface conditions \eqref{eq:1B} and \eqref{eq:1C}. Similar to the description in Section \ref{subsection: linear pde}, $\mathcal{L}(\vec{\alpha}_1, \vec{\alpha}_2)$ can be rearranged as 
\begin{equation*}
\text{Loss}(\vec{\alpha}_1,\vec{\alpha}_2)= \|\vec{\rho}\mathbf{A} \vec{\alpha} - \vec{\rho} \mathbf{b}\|_2^2.
\end{equation*}
Specifically, 
\begin{equation*}
\vec{\alpha} = \left[
 \begin{array}{cc}\vec{\alpha}_1\\ \vec{\alpha}_2\end{array}\right],\quad 
 \mathbf{A} = 
 \left[
 \begin{array}{cc}
  \mathbf{A}_{11}, &\mathbf{0}\\
  \mathbf{0}, &\mathbf{A}_{22} \\ 
 \mathbf{0}, &\mathbf{A}_{32} \\ 
  -\mathbf{A}_{41} , &\mathbf{A}_{42} \\ 
 -\mathbf{A}_{51} , &\mathbf{A}_{52} \\ 
 \end{array}
\right],  \quad \mathbf{b} =  \left[\begin{array}{cc}\mathbf{b}_1\\ \mathbf{b}_2\\ \mathbf{b}_3\\ \mathbf{b}_4 \\ \mathbf{b}_5  
\end{array}\right],
\end{equation*}
where
\begin{equation*}
\mathbf{A}_{11}= \left[
 \begin{array}{ccc}
  \mathcal{L}_1\left( t^1_1(\mathbf{x}_1^{d_1})\right) &\cdots & \mathcal{L}_1\left( t^1_I(\mathbf{x}_1^{d_1})\right)\\
  \vdots & \ddots & \vdots\\
   \mathcal{L}_1\left( t^1_1(\mathbf{x}_{n_{d_1}}^{d_1})\right) &\cdots & \mathcal{L}_1\left( t^1_I(\mathbf{x}_{n_{d_1}}^{d_1})\right)
 \end{array}
\right],
\
\
\mathbf{A}_{22}= \left[
 \begin{array}{ccc}
  \mathcal{L}_2\left( t^2_1(\mathbf{x}_1^{d_2})\right) &\cdots & \mathcal{L}_2\left( t^2_I(\mathbf{x}_1^{d_2})\right)\\
  \vdots & \ddots & \vdots\\
   \mathcal{L}_2\left( t^2_1(\mathbf{x}_{n_{d_2}}^{d_2})\right) &\cdots & \mathcal{L}_2\left( t^2_I(\mathbf{x}_{n_{d_2}}^{d_2})\right)
 \end{array}
\right],
\end{equation*}

\begin{equation*}
 \mathbf{A}_{32} = \left[
 \begin{array}{ccc}
  t^2_1(\mathbf{x}_1^{b})  &\cdots &   t^2_I(\mathbf{x}_1^{b}) \\
  \vdots & \ddots & \vdots\\
   t^2_1(\mathbf{x}_{n_{b}}^{b})  &\cdots &  t^2_I(\mathbf{x}_{n_{b}}^{b}) 
 \end{array}
\right],\
 \mathbf{A}_{41} = \left[
 \begin{array}{ccc}
 t^1_1(\mathbf{x}_1^{\gamma}) &\cdots & t^1_I(\mathbf{x}_1^{\gamma}) \\
  \vdots & \ddots & \vdots\\
  t^1_1(\mathbf{x}_{n_\gamma}^{\gamma}) &\cdots & t^1_I(\mathbf{x}_{n_\gamma}^{\gamma})
 \end{array}
\right],
\
 \mathbf{A}_{42} = \left[
 \begin{array}{ccc}
 t^2_1(\mathbf{x}_1^{\gamma}) &\cdots & t^2_I(\mathbf{x}_1^{\gamma}) \\
  \vdots & \ddots & \vdots\\
  t^2_1(\mathbf{x}_{n_\gamma}^{\gamma}) &\cdots & t^2_I(\mathbf{x}_{n_\gamma}^{\gamma})
 \end{array}
\right],
\end{equation*}
\begin{equation*}
  \mathbf{A}_{51} = \left[
 \begin{array}{ccc}
 a_1 \nabla t^1_1(\mathbf{x}_1^{\gamma})\cdot \mathbf{n} &\cdots & a_1 \nabla 
 t^1_I(\mathbf{x}_1^{\gamma})\cdot \mathbf{n} \\
  \vdots & \ddots & \vdots\\
  a_1 \nabla t^1_1(\mathbf{x}_{n_\gamma}^{\gamma})\cdot \mathbf{n} &\cdots & a_1 \nabla 
 t^1_I(\mathbf{x}_{n_\gamma}^{\gamma})\cdot \mathbf{n}
 \end{array}
\right],
\
 \mathbf{A}_{52} = \left[
 \begin{array}{ccc}
 a_2 \nabla t^2_1(\mathbf{x}_1^{\gamma})\cdot \mathbf{n} &\cdots & a_2 \nabla 
 t^2_I(\mathbf{x}_1^{\gamma})\cdot \mathbf{n} \\
  \vdots & \ddots & \vdots\\
  a_2 \nabla t^2_1(\mathbf{x}_{n_\gamma}^{\gamma})\cdot \mathbf{n} &\cdots & a_2 \nabla 
 t^2_I(\mathbf{x}_{n_\gamma}^{\gamma})\cdot \mathbf{n}
 \end{array}
\right],
\end{equation*}
\begin{equation*} 
\mathbf{b}_1 = \left[
 \begin{array}{c}
  f_1(\mathbf{x}_1^{d_1})\\
  \vdots \\
  f_1(\mathbf{x}_{n_{d_1}}^{d_1}) \\
 \end{array}
\right],\
\mathbf{b}_2 = \left[
  \begin{array}{c}
  f_2(\mathbf{x}_1^{d_2})\\
  \vdots \\
  f_2(\mathbf{x}_{n_{d_2}}^{d_2}) \\
 \end{array}
\right],\
\mathbf{b}_3 = \left[
  \begin{array}{c}
  h(\mathbf{x}_1^{b})\\
  \vdots \\
  h(\mathbf{x}_{n_{b}}^{b}) \\
 \end{array}
\right]
\
\mathbf{b}_4 = \left[
  \begin{array}{c}
  g_d(\mathbf{x}_1^{\gamma})\\
  \vdots \\
  g_d(\mathbf{x}_{n_{\gamma}}^{\gamma}) \\
 \end{array}
\right],
\
\mathbf{b}_5 = \left[
  \begin{array}{c}
  g_n(\mathbf{x}_1^{\gamma})\\
  \vdots \\
  g_n(\mathbf{x}_{n_{\gamma}}^{\gamma}) \\
 \end{array}
\right],
\end{equation*}
and $\vec{\rho}= \texttt{diag}\left[\sqrt{\frac{\rho_1^{d_1}}{n_{d_1}}}, \cdots, \sqrt{\frac{\rho_{n_{d_1}}^{d_1}}{n_{d_1}}},\sqrt{\frac{\rho_1^{d_2}}{n_{d_2}}}, \cdots, \sqrt{\frac{\rho_{n_{d_2}}^{d_2}}{n_{d_2}}},\sqrt{\frac{\rho_1^b}{n_{b}}}, \cdots, \sqrt{\frac{\rho_{n_{b}}^b}{n_{b}}},\sqrt{\frac{\rho_1^{\gamma_d}}{n_{\gamma}}}, \cdots, \sqrt{\frac{\rho_{n_\gamma}^{\gamma_d}}{n_\gamma}},\sqrt{\frac{\rho_1^{\gamma_n}}{n_{\gamma}}}, \cdots, \sqrt{\frac{\rho_{n_\gamma}^{\gamma_n}}{n_\gamma}}\right]$.
The specific penalty weights are set as $\rho_j^{d_k} = \frac{n_{d_k}}{ \sum_{i=1}^I \mathcal{L}_k(t_i^k( \mathbf{x}_j^{d_k}))^2}$ with $1\leq j\leq n_{d_k}$ and $k=1,2$, $\rho_j^{b} = \frac{n_{b}}{ \sum_{i=1}^I (t_i^2( \mathbf{x}_j^{b}))^2}$ with $1\leq j\leq n_{b}$, and
$\rho_j^{\gamma_d} = \frac{n_{\gamma}}{ \sum_{i=1}^I (t_i^1( \mathbf{x}_j^{\gamma}))^2 +  \sum_{i=1}^I  (t_i^2( \mathbf{x}_j^{\gamma}))^2}, \ \rho_j^{\gamma_n} = \frac{n_{\gamma}}{ \sum_{i=1}^I (a_1 \nabla t^1_i(\mathbf{x}_j^{\gamma})\cdot \mathbf{n})^2 +  \sum_{i=1}^I  (a_2 \nabla t^2_i(\mathbf{x}_j^{\gamma})\cdot \mathbf{n})^2}$ with $1\leq j\leq n_{\gamma}$.

\section{Details of the training and test points setting\label{appendix:2}}
We specify the setup for Example 1 to Example 4 as follows:  
\begin{itemize}
    \item[-] In Examples \ref{example 1} and  \ref{example 2}: We evaluate the loss function \eqref{eq: loss function of fto-pinn} at 8000 randomly selected points, which are divided into three sets for evaluating the PDE conditions (4000), the Dirichlet boundary condition (2000), and the initial boundary condition (2000), respectively. The test errors are evaluated on a $129\times 129$ uniform mesh in $\Omega=[0, 1]^2$.

    \item[-] Example  \ref{example 3}: We evaluate the loss function \eqref{eq: loss function of fto-pinn} at 3101 randomly selected points, which are divided into three sets for evaluating the PDE conditions (2000), the Dirichlet boundary condition (1000), and the initial boundary condition (101). The test errors are evaluated on a $101\times 101$ uniform mesh in $\Omega=[0, 1]^2$.

    \item[-] Example  \ref{example 4}: We evaluate the loss function \eqref{eq: loss function of fto-pinn} at 4,000 randomly selected points, which are divided into three sets for evaluating the PDE conditions (2000), the Dirichlet boundary condition (1000), and the interface boundary condition (1000). The test errors are evaluated on a $101\times 101$ uniform mesh in $\Omega=[-1, 1]^2$.
   
\end{itemize}

\end{document}